\documentclass[11pt, psamsfonts]{amsart}
\usepackage{amsmath}
\usepackage{amsthm}
\usepackage{amssymb}
\usepackage{amscd}
\usepackage{amsfonts}
\usepackage{amsbsy}
\usepackage{epsfig}
\usepackage{graphicx}

\newtheorem{theorem}{Theorem}
\newtheorem{remark}{Remark}
\newtheorem{proposition}{Proposition}
\newtheorem{lemma}{Lemma}
\newtheorem{corollary}{Corollary}

\newenvironment{definition}{{\bf Definition.}}{}

\newfont\bbf{msbm10 at 12pt}
\def\eps{\varepsilon}
\def\phi{\varphi}
\def\R{{\mathbb R}}
\def\C{{\mathbb C}}
\def\N{{\mathbb N}}
\def\Z{{\mathbb Z}}

\def\J{{\mathcal J}}
\def\B{{\mathcal B}}
\def\I{{\mathcal I}}
\def\P{{\mathcal P}}
\def\F{{\mathcal F}}
\def\D{{\mathcal D}}

\def\cont{{\mathcal C}}
\def\es{{\emptyset}}
\def\sm{\setminus}

\def\lev{\mbox{level}}
\def\ulambda{\underline\lambda}

\def\tag#1{\hfill \qquad  #1}
\def\dist{\mbox{dist}}

\def\supp{\mbox{\rm supp}}
\def\modulus{\mbox{mod}}

\def\diam{\mbox{\rm diam} }
\def\orb{\mbox{\rm orb}}
\def\Crit{{\mathcal Cr}}
\def\cut{{\mathcal Cut}}

\def\bd{\partial }
\def\theta{\vartheta}
\def\le{\leqslant}
\def\ge{\geqslant}

\newcommand{\st}{such that }

\begin{document}
\bibliographystyle{plain}
\title[Markov Extensions and Measures for Complex Polynomials.]
{Markov Extensions and Lifting Measures for Complex Polynomials.}
\author{Henk Bruin, Mike Todd}
\thanks{This research was supported by EPSRC grant GR/S91147/01}

\begin{abstract} For polynomials $f$ on the
complex plane with a dendrite Julia set we study invariant
probability measures, obtained from a reference measure. To do
this we follow Keller \cite{keller} in constructing canonical
Markov extensions. We discuss `liftability' of measures (both
$f$--invariant and non--invariant) to the Markov extension,
showing that invariant measures are liftable if and only if they
have a positive Lyapunov exponent. We also show that
$\delta$--conformal measure is liftable if and only if the set of
points with positive Lyapunov exponent has positive measure.
\end{abstract}

\maketitle

\section{Introduction}\label{sec_intro}

Ergodic properties for polynomial or rational maps have been
looked at for various measures and various types of Julia sets.
One can consider the measures of maximal entropy, e.g.
\cite{FLM,Zdu}, or more generally, equilibrium states of certain
H\"older potentials, see for example \cite{lyrational, DPU,
haydn}. This approach is particularly natural when the map is
hyperbolic and the potential is $-t\log |Df|$, where $t$ is the
Hausdorff dimension of the Julia set: then the equilibrium state
is equivalent to conformal measure (as obtained by Sullivan, see
\cite{sull}). When the Julia set is parabolic, invariant measures
equivalent to conformal measure are found in \cite{DUparabolic,
Utame, Unonrecurrent}. In the case where there are recurrent
critical points in the Julia set, the papers \cite{GS1,Prz,Rees}
focus on invariant probability measures that are absolutely
continuous with respect to conformal measures, using assumptions
on the derivatives on the critical orbits. See \cite{PU,Ureview}
for surveys.  The theory has not yet developed to the same extent
as, for example, interval maps, where the availability of induced
maps and tower constructions (cf. Young \cite{Young1,Young2})
allowed the investigation of several stochastic properties,
including the rate of mixing and Central Limit Theorem
\cite{Young2, BLS}, return time statistics and related properties
\cite{BSTV,BrV,Col} and Invariance Principles, see e.g. \cite{MN}.
However, during the preparation of this paper, we learned that
some good results
in this direction have been proved for rational maps satisfying
the `Topological Collet--Eckmann' condition in \cite{PrR-L}.

In the 1980s, Hofbauer and Keller constructed so--called canonical
Markov extensions for piecewise monotone maps of the interval
\cite{H1,HK,keller}, which they used to study the topological and
measure theoretical behaviour of these maps.  These Markov
extensions were considered in an abstract setting in
\cite{keller,Bumulti}, where one of the aims was to extend the
theory to higher dimensions. Indeed in \cite{Bumulti} some higher
dimensional examples are given. That paper focuses on the
probability measures given by the symbolic dynamics obtained from
the tower structure, an approach also used in \cite{New, BuSa}. In
\cite{BuK} results on transfer operators are proved in this higher
dimensional setting and in \cite{BuPS} conformal measures are
found.

Our approach follows the papers of Keller,
\cite{keller,Keller_Hopf}. In the first of these papers, results
are proved about the liftability of probability measures on the
original system to the associated Markov extension.  In
particular, the liftability of ergodic invariant measures with
positive entropy is shown. While the abstract theory given there
applies, in principal, in any dimension, the applications given
are to interval maps. In the second paper it is shown
that, given a smooth interval map, positive pointwise Lyapunov
exponents implies the liftability of Lebesgue measure.  The
purpose of this paper is to extend those results to maps on the
complex plane. We construct Markov extensions $(\hat \J, \hat f)$
for complex polynomials $f$ and study the liftability properties
of probability measures supported on the Julia set $\J$. (This
allows us to deal with some cases where critical points lie in
$\J$.) Given a probability measure $\mu$ on $\J$, we construct a
sequence of Cesaro means $\{ \hat \mu_n\}_n$ on $\hat \J$, and we
say that $\mu$ is {\em liftable} to the Markov extension if this
sequence has a non--zero vague limit measure $\hat \mu$.  The
limit measure $\hat \mu$ is $\hat f$--invariant, even if the
measure $\mu$ is not $f$--invariant. This technique is
particularly useful for finding invariant probability measures
that are absolutely continuous with respect to $\delta$--conformal
measure on the Julia set.

Among other things, we prove that an ergodic invariant probability
measure $\mu$ is liftable if and only if its Lyapunov exponent is
positive (cf. \cite{BK}). Furthermore, for liftable measures,
typical points are conical (i.e. go to large scale, see
Lemma~\ref{lem:large_scale}) with positive frequency. Similar
results hold for (non--invariant) $\delta$--conformal
measure $\mu_\delta$.
(The measure $\mu_\delta$ is {\em $\delta$--conformal}
on $\J$ if $\mu_\delta(\J)=1$ and $\mu(f(A)) = \int_A
|Df|^\delta \ d\mu_\delta$ for all measurable sets $A$ such that
$f:A\to f(A)$ is $1$--to--$1$.)
We prove the pointwise lower Lyapunov exponent
$\ulambda(z)$ is strictly greater than 0 for a set of positive
$\mu_\delta$--measure if and only if $\mu_{\delta}$ is liftable,
and in this case there is an $f$--invariant probability measure
equivalent to $\mu_\delta$.  We note that this result applies to
polynomials considered in \cite{GS1,Prz,Rees}, when the Julia sets
of these polynomials are dendrites, see below.

When proving our results on the relation between liftability and
positive Lyapunov exponents, we use the Koebe Lemma: a
one--dimensional tool.  Work in progress aims at extending these
results to higher dimensions. Our result on finding invariant
probability measure absolutely continuous with respect to
$\delta$-conformal measure again uses the Koebe Lemma and seems a
more difficult type of result to generalise.

Markov extensions (popularly called Hofbauer towers) are less well
known than the Young towers, \cite{Young1,Young2}.  We wish to
highlight the difference between these two constructions.  In
short, for the Young tower case, given an invariant measure $\mu$
and a subset $Y$ of the phase space, a partition $Y = \cup_j Y_j
\pmod \mu$ is constructed together with return times $R_j$ such
that $F:\cup_j Y_j \to Y$, $F|_{Y_j} = f^{R_j}|_{Y_j}$ and $f^{R_j}:Y_j
\to Y$ is $1$--to--$1$ and has good distortion and expansion
properties.  These are then used to study stochastic limit
properties (e.g. mixing rates, the Central Limit Theorem,
invariance principles) of specific invariant measures.  The
construction is therefore linked to the choice of the measure, and
may be quite involved in practical applications. The construction
of the Hofbauer tower, on the other hand, is combinatorial and can
be used to study all probability measures. In fact, it is exactly
for the liftable invariant probability measures that Young towers
can be constructed, in a canonical way, as first return maps to
appropriate sets in the Markov extension, see \cite{bruin1}.

The structure of this paper is as follows. The construction of the
canonical Markov extension occupies Section~\ref{sec_extension}.
We restrict our attention to polynomials $f$ with locally
connected full Julia sets (dendrites), as we need to find a finite
partition $\P_1$ of the Julia set $\J$ such that $f$ is univalent
on each partition element. Such partitions may exist for the Julia
set of many other rational maps as well, but is hard to give for
rational maps in all generality. In Section~\ref{sec_lifting} we
describe the lifting procedure of measures.  As remarked there, in
contrast to subsequent sections, Section~\ref{sec_lifting} is
largely independent of the geometry of $\J$, and can be easily
extended to Markov extensions in other settings.   In
Section~\ref{sec_induce} we introduce inducing constructions as a
tool to prove that for liftable measures, typical points will `go
to large scale' with positive frequency. Section~\ref{sec_lyap}
focuses on (ergodic) invariant probability measures $\mu$ and
their Lyapunov exponents $\lambda(\mu)$. It is shown that $\mu$ is
liftable if and only if $\lambda(\mu) > 0$.
Section~\ref{sec_conformal} gives a similar result for
$\delta$--conformal measure $\mu_\delta$. It shows that
$\mu_\delta$ is liftable if and only if the pointwise lower
Lyapunov exponent $\ulambda(z) > 0$ for all $z$ in a set of
positive $\mu_\delta$--measure.

{\bf Acknowledgements:} We are grateful to J.\ Hawkins, M.\
Urba\'nski and J.\ Rivera--Letelier for fruitful discussions. Also
the referee's suggestions for improving clarity of the paper are
gratefully acknowledged.

\section{The Markov Extension}\label{sec_extension}

Let $f:\C \to \C$ be a polynomial of degree $d$ with a connected,
locally connected and full Julia set $\J$ (i.e. $\C\setminus \J$
is connected). Consequently all critical points belong to $\J$.
Let $\Crit$ denote the critical set. It is easy to see that $\J$
is a dendrite, defined as follows (cf. \cite{kura}).

\begin{definition}
A metric space $(X,d)$ is called a \emph{dendrite} if it is
connected, locally connected, and for any two points $x,y\in X$
there is a unique arc $\gamma:[0,1]\to X$ connecting $x$ to $y$.
\end{definition}

The Fatou set $\F$ coincides with the basin of $\infty$. Let the
Green function $G:\F \to \R$ be defined by $G(z) = \lim_{n \to
\infty} \frac{\log |f^n(z)|}{d^n}$, see \cite{Mil} for more
details. The equipotentials (i.e. level sets) of the Green
function form a foliation of $\F$ consisting of nested Jordan
curves. The orthogonal foliation is the foliation of external
rays. Each external ray is a copy of $\R$ embedded in $\F$, and if
$\gamma:\R \to R$ is such an embedding such that $|\gamma(t)|$ is
large for large $t$, then $\lim_{t \to \infty} \arg \gamma(t)$ is
a well--defined number $\theta \in S^1$, called the external angle
of $R$. Let $R_{\theta}$ denote the ray with external angle
$\theta$, and $\gamma_{\theta}:\R \to R_{\theta}$ its
parameterisation. It is convenient to parameterise external rays
by the values of the Green function: $G(\gamma_{\theta}(t)) = t$
for each $\theta \in S^1$ and $t \in \R$. Note that $f(R_{\theta})
= R_{d \theta \bmod 1}$; more precisely: $f(\gamma_{\theta}(t)) =
\gamma_{d \theta \bmod 1} (t+1)$.

\begin{lemma}\label{Lempartition}
There is a finite partition $\P_1$ of $\J\sm\Crit$ such that
$f|_Z$ is univalent for each $Z \in \P_1$.
\end{lemma}

\begin{proof}
Because $\J$ is locally connected, each external ray $R_{\theta}$
lands at a single point in $z \in \J$ and each $z \in \J$ is the
landing point of at least one external ray. For each $c \in
\Crit$, select $\kappa_c \ge 1$ rays that land at $f(c)$ (note
that Theorems 1.1 and 3.1 of \cite{kiwi} imply that there can be
at most $2^d$ rays landing here). If $c$ has degree $d_c$, there
are $\kappa_c d_c \ge 2$ preimage rays landing at $c$. The union
of these rays, together with $c_j$, is a locally compact set,
separating the plane and also $\J$ into $\kappa_c d_c$ `segments'.
As $\J$ is closed, the closure of the segments of $\J \setminus \{
c \}$ intersect only at $c$. Repeating the argument for all
other critical points gives the assertion.
\end{proof}

Let $\P_0 = \{ \J \}$ be the trivial
partition of $\J$, and $\P_1$ be the partition of
Lemma~\ref{Lempartition}. Let $\P_n = \bigvee_{i = 0}^{n-1}
f^{-i}(\P_1)$, and for $z \notin \cup_{i=0}^{n-1} f^{-i}(\Crit)$ let
$Z_n[z]$ be the element of $\P_n$ containing $z$. Note that each
$Z_n\in \P_n$ is connected (in fact, $\overline{Z_n}$ is a
dendrite) and $f^n|_{Z_n}$ is univalent. We call $Z_n$ an
$n$--cylinder.

\begin{definition} The canonical Markov extension $\hat \J$ is a disjoint union
of copies $D$ of subsets of $\J$, subject to an identification
discussed below. We call sets $D$ {\em domains} and denote their
collection as by $\D$. Let $\pi:\hat \J \to \J$ be the inclusion
map. Domains $D$ are defined recursively as follows:
\begin{itemize}
\item The first domain $\hat \J_0$, called the {\em base} of the Markov
extension, is a copy of $\J$.

\item Given a domain $D\in \D$ and a nonempty
set of the form $\overline{ f(\pi(D) \cap Z)}$ for some  $Z \in
\P_1$, we let $D'$ be a copy of  $\overline{ f(\pi(D) \cap Z)}$
and add it to $\D$. Write $D \to D'$ in this case.
\item The collection $\D$ is such that $\hat \J_0 \in \D$, and $\D$ is closed
under the previous operation.
\end{itemize}
Each $\hat z\in D$ can be represented by a pair $(z,D)$ where
$\hat z\in D$ and $\pi(\hat z)=z$.  Moreover, any pair $(z,D)$
defines a unique $\hat z\in \J$ whenever $z\in \pi(D)$.  This
allows us to define $\hat f:\hat \J \to \hat \J$.
\begin{itemize}
\item If $\hat z \in D$, $D \to D'$ and $\pi(\hat z)$ belongs to the
closure of $Z \in \P_1$ such that $\pi(D') = \overline{ f(\pi(D)
\cap Z)}$, then we let $\hat f(\hat z) = (f(z), D')$. Clearly
\[
\pi \circ \hat f = f \circ \pi.
\]
If $\pi(\hat z)\in \Crit$, then $\hat f$ can be multi-valued at
$\hat z$, but a domain $D\in \D$ contains at most one of the
images of $\hat z$. In all other cases, $\hat f(\hat z)$ is a
single point, belonging to a single domain $D$.
\end{itemize}
The next step is to define the {\em cutpoints}, their {\em ages} and
{\em origins},
as well as the {\em level} of domains.
\begin{itemize}
\item  The base $\hat \J_0$ contains no cutpoints.
\item If $\hat z \in D$ is a cutpoint or $\pi(\hat z) \in \Crit$,
then each image $\hat f(\hat z)$ is a cutpoint. Its \emph{age} is
\[
\left\{ \begin{array}{cl}
1 & \mbox{ if } \pi(\hat z) \in \Crit \mbox{ and } \hat z \mbox{ is not a cutpoint;}\\
a+1 & \mbox{ if } \hat z \mbox{ is a cutpoint of age } $a$.
\end{array} \right.
\]
The set of cutpoints is denoted by $\cut$.
\item An \emph{$a$--cutpoint} will be a cutpoint of age $a$.
Each $a$--cutpoint $\hat z$ satisfies $\hat z = (f^a(c),D)$ for
some $D \in \D$ and $c \in \Crit$. This critical point $c$ is
called the {\em origin} of $\hat z$.
\item  Given a domain $D$, $\lev(D)$ is $0$ if there are no
cutpoints in $D$, and is the maximal age of the cutpoints in $D$
otherwise. Let $\hat \J_R$ be the union of all domains of $\lev(D)
\le R$.
\end{itemize}
The final step is the identification of domains wherever possible.
\begin{itemize}
\item
Any two domains $D$ and $D'$ such that $\pi(D) = \pi(D')$,
$\pi(D\cap\cut)=\pi(D'\cap\cut)$ and  whose cutpoints have the
same ages and origins are identified. The {\em canonical Markov
extension} is the disjoint union of the domains, factorised over
the identification described above.
\end{itemize}
The arrow relations $D \to D'$ give the Markov extension the
structure of an (infinite) directed graph. This is the {\em Markov
graph}, since by construction $(\hat \J,  \hat f)$ is Markov with
respect to the partition of the domains of $\hat \J$.
\end{definition}

For counting arguments later (see Lemma~\ref{lem:basiccounting}
and the appendix), we must be aware of the possibility of `moving
sideways' in the Markov graph. That is, it is possible that for
some domain $D$ of $\hat \J$ there is an arrow $D \to D'$ where
$\lev(D) = \lev(D')$.  This occurs if $D$ is a domain of level $n$
containing one cutpoint $\hat z$ of age $n$ and a cutpoint $\hat
z'$ of age $n-1$. If the arc in $D$ connecting these two cutpoints
intersects $\pi^{-1}(\Crit)$ (recall that since $\J$ is a
dendrite, for any $z,z'\in \J$ there exists a unique arc in $\J$
connecting $z$ to $z'$), then the domain $D'$ containing $\hat
f(\hat z)$ will also have level $n$. So if $D \to D'$, then
$\lev(D')$ can take \emph{any} value $\le \lev(D)+1$.

Define $\hat\P_1$ to be the partition given by $\D \vee\pi^{-1}
\P_1$.
Let $\hat \P_n:= \bigvee_{i=0}^{n-1}\hat f^{-i}(\P_1)$ and
$\hat\P_n^R:= \hat \P_n \cap \hat \J_R$.

\begin{remark}
The partition in Lemma~\ref{Lempartition}, and hence the
construction of the Markov extension, is not unique,
because we have freedom in choosing the number of the rays
$\kappa_c$ for each critical value. However, any choice makes a
valid partition. To illustrate this,  assume that $f(z) = z^2+c$
for $c \in (-2,-\frac14)$ such that $0 \in \J$. One is inclined to
choose two (complex conjugate) rays landing at $c = f(0)$, see
Figure~\ref{FigRays} (left). This will lead to a canonical Markov
extension which is very similar to the standard Markov extension
constructed for interval maps. More precisely, select the domains
$D \in \D$ such that $\pi(D) \cap \R \neq \es$ and such that
if $\pi(D) = f^n(\overline{Z_n})$, then for each $x \in \pi(D)
\cap \R$, there is $x_0 \in \overline{Z_n} \cap \R$ such that $x =
f^n(x_0)$. For each such $D$, retain $D \cap \pi^{-1}(\R)$, and
discard the rest of $D$ as well as all other domains. Then this
set with remaining graph structure is exactly the real Markov
extension, see Figure~\ref{FigRays} (left, bold lines).
\par
Choosing only one ray is possible as well; in this case, each
domain in the Markov extension will be a copy of the whole Julia
set, see Figure~\ref{FigRays} (right), and they will be
distinguished only by the fact that they have different (numbers of) cutpoints,
and consequently different canonical neighbourhoods, see below.
\end{remark}

\begin{figure}[ht]
\begin{center}
\begin{minipage}{140mm}
\unitlength=6mm
\begin{picture}(22,22)(-2,-3)
\let\ts\textstyle
\thinlines \put(10,-2.5){\line(0,1){21}}
\thicklines
\put(1,1){\line(1,0){8}}\put(1,0.95){\line(1,0){8}}\put(1,1.03){\line(1,0){8}}
\put(5,-1){\line(0,1){4}} \put(3,0){\line(0,1){2}}
\put(7,0){\line(0,1){2}} \put(2,0.5){\line(0,1){1}}
\put(4,0.5){\line(0,1){1}} \put(6,0.5){\line(0,1){1}}
\put(8,0.5){\line(0,1){1}} \put(4.5,2){\line(1,0){1}}
\put(4.5,0){\line(1,0){1}} \thinlines \put(0, 0){$\hat \J_0$}
\put(2.5,1){\circle*{0.1}} \put(2.5, 1.2){\scriptsize $\hat c$}
\put(1,1){\circle*{0.1}} \put(0.7, 1.2){\scriptsize $-\beta$}
\put(9,1){\circle*{0.1}} \put(8.9, 1.2){\scriptsize $\beta$}
\put(2.5,1){\line(-1,2){1.5}}\put(0.3, 4){\scriptsize $R_\phi$}
\put(2.5,1){\line(-1,-2){1.5}}\put(0.3, -2){\scriptsize
$R_{\phi'}$} \put(5,1){\line(-1,1){3}}\put(2, 4){\scriptsize
$R_{\phi'/2}$} \put(5,1){\line(-1,-1){3}}\put(8, 4){\scriptsize
$R_{\phi/2}$} \put(5,1){\line(1,1){3}}\put(2.4, -2){\scriptsize
$R_{(1+\phi)/2}$} \put(5,1){\line(1,-1){3}}\put(5.6,
-2){\scriptsize $R_{(1+\phi')/2}$}

\thicklines \put(11,1){\line(1,0){8}} \put(15,-1){\line(0,1){4}}
\put(13,0){\line(0,1){2}} \put(17,0){\line(0,1){2}}
\put(12,0.5){\line(0,1){1}} \put(14,0.5){\line(0,1){1}}
\put(16,0.5){\line(0,1){1}} \put(18,0.5){\line(0,1){1}}
\put(14.5,2){\line(1,0){1}} \put(14.5,0){\line(1,0){1}} \thinlines
\put(19, 0){$\hat \J_0$} \put(12.5,1){\circle*{0.1}} \put(12.5,
1.2){\scriptsize $\hat c$} \put(11,1){\circle*{0.1}} \put(10.7,
1.2){\scriptsize $-\beta$} \put(19,1){\circle*{0.1}} \put(18.9,
1.2){\scriptsize $\beta$} \put(12.5,1){\line(-1,2){1.5}}\put(10.3,
4){\scriptsize $R_\phi$} \put(15,1){\line(-1,-1){3}}\put(18,
4){\scriptsize $R_{\phi/2}$} \put(15,1){\line(1,1){3}}\put(12.4,
-2){\scriptsize $R_{(1+\phi)/2}$}

\thicklines
\put(2.5,8){\line(1,0){6.5}}\put(2.5,7.95){\line(1,0){6.5}}\put(2.5,8.03){\line(1,0){6.5}}
\put(0.5,8){\line(1,0){1.5}} \put(5,6){\line(0,1){4}}
\put(3,7){\line(0,1){2}} \put(7,7){\line(0,1){2}}
\put(1.5,7.5){\line(0,1){1}} \put(4,7.5){\line(0,1){1}}
\put(6,7.5){\line(0,1){1}} \put(8,7.5){\line(0,1){1}}
\put(4.5,9){\line(1,0){1}} \put(4.5,7){\line(1,0){1}} \thinlines
\put(2.5,8){\circle*{0.3}} \put(2.5, 8.2){\scriptsize $\hat c$}
\put(2,8){\circle*{0.3}} \put(1.95, 8.2){\scriptsize $\hat c$}
\put(0.5,8){\circle*{0.1}} \put(0.2, 8.2){\scriptsize $-\beta$}
\put(9,8){\circle*{0.1}} \put(8.9, 8.2){\scriptsize $\beta$}
\put(2.5,8){\line(-1,2){1.5}}\put(1, 11){\scriptsize $R_\phi$}
\put(2.5,8){\line(-1,-2){1.5}}\put(1.1, 5){\scriptsize
$R_{\phi'}$} \put(2,8){\line(-1,2){1.5}}\put(-0.2, 11){\scriptsize
$R_\phi$} \put(2,8){\line(-1,-2){1.5}}\put(-0.3, 5){\scriptsize
$R_{\phi'}$}

\thicklines \put(11,8){\line(1,0){8}} \put(15,6){\line(0,1){4}}
\put(13,7){\line(0,1){2}} \put(17,7){\line(0,1){2}}
\put(12,7.5){\line(0,1){1}} \put(14,7.5){\line(0,1){1}}
\put(16,7.5){\line(0,1){1}} \put(18,7.5){\line(0,1){1}}
\put(14.5,9){\line(1,0){1}} \put(14.5,7){\line(1,0){1}} \thinlines
\put(12.5,8){\circle*{0.3}} \put(12.5, 8.2){\scriptsize $\hat c$}
\put(11,8){\circle*{0.1}} \put(10.2, 8.2){\scriptsize $-\beta$}
\put(19,8){\circle*{0.1}} \put(18.9, 8.2){\scriptsize $\beta$}
\put(12.5,8){\line(-1,2){1.5}}\put(11, 11){\scriptsize $R_\phi$}

\thicklines
\put(1.5,15){\line(1,0){4}}\put(1.5,14.95){\line(1,0){4}}\put(1.5,15.03){\line(1,0){4}}
\put(2,14){\line(0,1){2}} \put(6.5,15){\line(1,0){2.5}}
\put(4,13){\line(0,1){4}} \put(7,14){\line(0,1){2}}
\put(3,14.5){\line(0,1){1}} \put(5,14.5){\line(0,1){1}}
\put(8,14.5){\line(0,1){1}} \put(3.5,16){\line(1,0){1}}
\put(3.5,14){\line(1,0){1}} \thinlines \put(1.5,15){\circle*{0.3}}
\put(1.5, 15.2){\scriptsize $\hat c$} \put(5.5,15){\circle*{0.3}}
\put(5.6, 15.2){\scriptsize $\hat f(\hat c)$}
\put(6.5,15){\circle*{0.3}} 
\put(9,15){\circle*{0.1}} \put(8.9, 15.2){\scriptsize $\beta$}
\put(1.5,15){\line(-1,2){1.5}}\put(0, 18){\scriptsize $R_\phi$}
\put(1.5,15){\line(-1,-2){1.5}}\put(0.1, 12){\scriptsize
$R_{\phi'}$} \put(5.5,15){\line(1,3){1}}\put(5.5, 18){\scriptsize
$R_{2\phi'}$} \put(5.5,15){\line(1,-3){1}}\put(5.3,
12){\scriptsize $R_{2\phi}$} \put(6.5,15){\line(1,3){1}}\put(7.5,
18){\scriptsize $R_{2\phi'}$}
\put(6.5,15){\line(1,-3){1}}\put(7.5, 11){\scriptsize $R_{2\phi}$}

\thicklines \put(11,15){\line(1,0){8}} \put(15,13){\line(0,1){4}}
\put(13,14){\line(0,1){2}} \put(17,14){\line(0,1){2}}
\put(12,14.5){\line(0,1){1}} \put(14,14.5){\line(0,1){1}}
\put(16,14.5){\line(0,1){1}} \put(18,14.5){\line(0,1){1}}
\put(14.5,16){\line(1,0){1}} \put(14.5,14){\line(1,0){1}}
\thinlines \put(12.5,15){\circle*{0.3}} \put(12.5,
15.2){\scriptsize $\hat c$} \put(16.5,15){\circle*{0.3}} \put(16,
15.2){\scriptsize $\hat f(\hat c)$} \put(11,15){\circle*{0.1}}
\put(10.2, 15.2){\scriptsize $-\beta$} \put(19,15){\circle*{0.1}}
\put(18.9, 15.2){\scriptsize $\beta$}
\put(12.5,15){\line(-1,2){1.5}}\put(11, 18){\scriptsize $R_\phi$}
\put(16.5,15){\line(1,-3){1}}\put(16.3, 12){\scriptsize
$R_{2\phi}$}
%
%
\put(7.5,2.5){\vector(-1,3){1.5}}
\put(4.9,3.4){\vector(-1,1){3.8}} \put(14,2.5){\vector(0,1){3.5}}
\put(1.3,8.7){\vector(1,1){6}}
\put(8,7){\line(0,-1){1}}\put(8,6){\line(-1,0){1.5}}
\put(6.5,6){\vector(0,1){1}} \put(4.7, 10){\line(-2,1){2}}
\put(2.7, 11){\vector(-2,-3){1.5}}
\put(3.5,9.5){\vector(0,1){3.5}} \put(14,9.5){\vector(0,1){3.5}}

\put(18,6){\line(-1,0){1.5}} \put(16.5,6){\vector(0,1){1}}
\put(18.0,6){\line(0,1){1}}

\end{picture}
\caption{\small Schematic picture of some domains of the Markov
extension of a quadratic map $z \mapsto z^2 + c$
using (left) two external rays
$R_\phi$ and $R_\phi'$ for $\phi' = 2\pi-\phi$, landing at $\hat
c$ and (right) one external ray $R_\phi$.
\newline
The pictures on the bottom line are $\hat \J_0$ and the rays
landing at $0$ defining the partition are shown. The point
$\pi(\beta) > 0$ is fixed under $f$, and $\pi(-\beta)$ is its
other preimage.
\newline
The middle line gives domains of level $1$ and
the top line domains of level $2$. Cutpoints are denoted by a
$\bullet$
\newline
Arrows indicate the edges $D \to D'$ of the Markov graph.
For clarity of the
picture, if $Z$ and $Z' \in \P_1$ are symmetric to each other, the
arrow from only one of them is shown.
\newline
The bold lines on the left pictures indicate the `real Markov
extension'. \label{FigRays}}
\end{minipage}
\vskip-20pt
\end{center}
\end{figure}

We summarise some properties of $\hat\J$ in the following lemma.

\begin{lemma}
\begin{enumerate}
\item[(a)] for any $a \ge 1$, each $D\in \D$ contains
at most $\#\Crit$ cutpoints of age $a$ (and at most one for each
different origin $c$);
\item[(b)] Let $D$ and $D'$ be domains in $\hat \J$ of the same level,
sharing a cutpoint of maximal age,
i.e., $\hat p \in D$ and $\hat p' \in D$ are cutpoints of age
$a = \lev(D) = \lev(D')$ and $\pi(\hat p) = \pi(\hat p')$.
Suppose also that $\hat p$ and $\hat p'$ have the same origin.
Then $\pi(D) = \pi(D')$ or  $\pi(D) \cap \pi(D') = \pi(\hat p)$;
\item[(c)] The number of domains of level $l$ is bounded
by $\#\Crit \prod_c \kappa_c$
Consequently, the number of domains in $\hat \J_R$
is at most $1+R \#\Crit \prod_c \kappa_c$.
\end{enumerate}
\label{lem:basiccounting}
\end{lemma}

Notice that (b) implies that within a given level we can only
`move sideways' a uniformly bounded number of times.

\begin{proof}
For the proof of (a), note that if $\hat p\in \cut \cap D$ has age $a$,
then $\pi(\hat p) \in f^a(\Crit)$. As $\hat \J$ contains no loops, only one
such point exists for each $c$ and $a$.
So there are at most $\# \Crit$ cutpoints of age $a$.

To prove (b), let $D$ and $D'$ be as in the statement, and assume that
$\pi(\hat p) = \pi(\hat p') =: p = f^a(c)$, where $a$ is the age of $p$
and $\hat p$ and $c$ their common origin.
This means that there are dendrites $E$ and $E'$ intersecting at $f(c)$
such that $\pi(D) = f^{a-1}(E)$, $\pi(D) = f^{a-1}(E')$,
and $f^{a-1}|_E$ and $f^{a-1}|_{E'}$ are homeomorphic.

Assume first that  $E$ and $E'$ have at least an arc in common.
If $E \neq E'$, say $x \in E \setminus E'$, then there is $y$ such that
$[x,f(c)] \cap E' = [y,f(c)]$.
Here $[a,b]$ indicates the unique arc in $E$ connecting $a$ and $b$.
By construction, each set $\overline{\J \setminus \pi(D')} \cap \pi(D')$
consists of post--critical points, and  the same holds for
$\overline{\pi(D) \setminus \pi(D')} \cap \pi(D')$.
Since $f^{a-1}(E') = \pi(D')$, we have $f^n(y) \cap \Crit \neq \es$ for
some $n \in \Z$. There are two possibilities:
\begin{itemize}
\item $y \in \cup_{n \ge 2} f^n(\Crit)$, but then $D'$ must
have a cutpoint of age $> a$, contradicting maximality of $a$.\\
\item  $y \in \cup_{n \le 1} f^n(\Crit)$.
In this case there is $\tilde c \in \Crit$ and $0 \le s < a$
such that $f^{s-1}(y) \owns \tilde c$.
Take $y$ such that $s$ is maximal with this property.
Now $f^s(E)$ and $f^s(E')$ belong to the same of the sectors defined by the
$\kappa_{\tilde c}$ external rays landing at $f(\tilde c)$.
But then $f^s(E)$ and $f^s(E')$ both contain $f^s(x)$, a contradiction.
Consequently, $\pi(D) = \pi(D')$.
\end{itemize}
Otherwise, $E$ and $E'$ intersect only at $f(c)$.
First assume that $\orb(c)$ contains no further critical points.
Then $f^{a-1}$ is locally univalent at $f(c)$.
Thus if $E \cap E' = f(c)$, then (as $\J$ is a dendrite, containing no loops)
$f^{a-1}(E) \cap f^{a-1}(E) = p$.
The final case is that there is $s$ and $\tilde c
\in \Crit$ such that $f^{s-1}(p) = \tilde c$ and $f^s(E) \cap f^s(E')$
contains more than just $f^s(p) = f(\tilde c)$.
Then the previous argument shows that $f^s(E) = f^s(E')$,
and we again obtain  $\pi(D) = \pi(D')$.

Now to prove (c), note that for each $c \in \Crit$ ad domain $D$,
$\pi^{-1}(c) \cap D$ has at most $\kappa_c$ images under $\hat f$.
Under further iteration of $\hat f$, this number does not
increase, unless $f^s(c) = \tilde c$ for some $\tilde c \in \Crit$,
in which case the number of images can multiply by at most $\kappa_{\tilde c}$.
The worst case is that there are $\prod_c \kappa_c$ images.
In other words, for each $l$, there can be at most $\prod_c \kappa_c$
domains $D$ of level $l$ for which $\pi(D)$ pairwise
intersect only at $f^l(c)$.
Using (b), this gives at most $\#\Crit \prod_c \kappa_c$ domains
of level $l$ altogether.
\end{proof}

\begin{lemma}\label{Lemfibersmeet}
Given $\hat z, \hat z' \in \pi^{-1}(z)$, one of the following
three cases occurs:
\begin{enumerate}
\item[(a)] There exists $n$ such
that $\hat f^n(\hat z) = \hat f^n(\hat z')$;
\item[(b)] $\cap_n Z_n[z]$ has
positive diameter;
\item[(c)] $z \in \bigcup_{c \in \Crit}\bigcup_{n \in
\Z} f^n(c)$.
\end{enumerate}
In the latter two cases at least one of $\hat z, \hat z'$ visits
any $\hat \J_R$ only finitely often.
\end{lemma}

\begin{proof}
Assume that $z$ is not precritical, and that neither $\hat z$ nor
$\hat z'$ is a cutpoint. Let $D$ and $D'$ be such that $\hat z \in
D$ and $\hat z' \in D'$. If there is $n$ such that the cylinder
$Z_n[z]$ is contained in $\pi(D)$ as well as in $\pi(D')$, then
$\hat f^n(D \cap \pi^{-1}(Z_n[z])) = \hat f^n(D' \cap
\pi^{-1}(Z_n[z]))$, and $\hat f^n(\hat z) = \hat f^n(\hat z')$ as
in case (a).

If on the other hand there is no such $n$, then $Z := \cap_m
Z_m[z]$ has positive diameter as in case (b). Furthermore, $Z$
contains no critical point in its interior (here we mean interior
with respect to the relative
topology on $\J$), and if $\pi(D) \not\supset Z$, then
$\hat Z := D \cap \pi^{-1} Z$ contains a cutpoint of $D$. Let $p$
be such a cutpoint of maximal age, say $a$. Then $\hat f^k(\hat Z)
\owns \hat f^k(p)$ which has age $a+k$. It follows that all the
sets $\hat f^k(\hat Z)$, $k \ge 0$ are disjoint. As a result $\hat
Z$ can remain in $\hat \J_R$ for at most $R$ iterates.

If $\hat z$ (or $\hat z'$) is a cutpoint of age $a$ then we are in
case (c) and the age of $\hat z$ will increase under iteration of
$\hat f$. So for any $R$, $\hat f^{R+k}(\hat z)$ is outside $\J_R$
for any $k \ge 1$.

Finally, if $z$ is precritical, then we are in case (c) again.  It
is possible that $z$ belongs to the common boundary of several
cylinder sets $Z_n$, and $\hat f^{n+1}$ is multivalued at $\hat z$
and $\hat z'$. But each image $\hat f^{n+1}(\hat z)$ and $\hat
f^{n+1}(\hat z')$ is a cutpoint of its level, so it will
eventually climb in the Markov extension.
\end{proof}

{\bf Running assumptions:} We will  repeatedly invoke the
following assumptions on measures $\mu$, to almost surely rule out
cases (b) and (c) of Lemma~\ref{Lemfibersmeet}, as explained in
the next section. Typical cylinders should shrink:
\begin{equation}\label{SC}
\diam \ Z_n[z] \to 0 \hbox{ as } n \to \infty \hbox{ for }
\mu\mbox{--a.e. } z \in \J. \tag{SC}
\end{equation}
and the mass on the precritical points is $0$:
\begin{equation}\label{CO}
\mu(\cup_{n \le 0} f^n(\Crit)) = 0.
\tag{$\Crit_0$}
\end{equation}

By Theorem 3.2 of \cite{blokh_levin}, \eqref{SC} automatically
follows from our assumption that $f$ is a polynomial and $\J$ is
locally connected and full.  However, as we believe Markov
extensions can be of use also when $\J$ is not locally connected,
(in which case one should think of a different partition $\P_1$
than the one based on external rays landing at $\Crit$), we will
refer to this property whenever we use it.

{\bf Canonical neighbourhoods:} let $U_{\hat \J_0}$ be a copy of
the neighbourhood $U_{\J}$ of $\J$ bounded by the equipotential
$\{ G(z) = 0\}$. This is the {\em canonical neighbourhood} of
$\hat \J_0$. We will define a canonical neighbourhood $U_D$ for
each $D \in \D$; they are copies of subsets of $\C$. The inclusion
map $\pi$ is extended to $U_D$ in the natural way. In the proof of
Lemma~\ref{Lempartition} we chose $\kappa_c$ external rays landing
at the critical value $f(c)$. The preimage rays landing at
critical points (together with $\Crit$) divide $\pi(U_{\hat
\J_0})$ into $\#\P_1$ regions. The closure of each such region $O$
contains exactly one element of $\P_1$: if $Z \in \P_1$, let $O_Z$
be the corresponding region. For each $D = \overline{f(Z)}$, $Z
\in \P_1$, let $U_D$ be a copy of $f(O_Z) \cap U_{\J}$. This set
is bounded by external rays landing at critical values and by the
equipotential $\{ G(z) = 0\}$. We call $U_D$ the {\em canonical
neighbourhood} of $D$, although it is not a neighbourhood in the
strict sense: $D \setminus U_D$ consists of the cutpoints of $D$.

We continue recursively. If $D \to D'$ and $U_D$ is the canonical
neighbourhood of $D$, then $U_{D'}$ is a copy of $f(\pi(U_D) \cap
O_Z) \cap U_{\J}$, where $Z \in \P_1$ is such that $\pi(D') =
\overline{f(\pi(D) \cap Z)}$. It is bounded by external rays
landing at cutpoints in $D'$ and by $\{ G(z) = 0\}$.

Let $\hat U$ be the disjoint union of all canonical
neighbourhoods. Then $\hat f$ naturally extends univalently to
$\hat U$ by
\[
\hat f(z,U_D) = (f(z), U_{D'})
\]
if $z \in \overline{\pi(U_D) \cap O_Z}$ where $Z$ is such that $\pi(D') =
\overline{ f(Z \cap D)}$.

\begin{lemma}
The recursive definition of $U_D$ is independent of the path $\hat
\J_0 \to \dots \to D$ by which $D$ is reached.
\end{lemma}

\begin{proof}
Since no path leads into $\hat \J_0$, its canonical neighbourhood
is uniquely defined. Now take $D \in \D$, $D \neq \hat \J_0$ with
at least two arrows leading to $D$. (To prove the lemma, we can
restrict to domains $D$ with two arrows rather than two paths
leading to it, because when two paths eventually merge, it
suffices to study those domains at which these paths merge.) For
any cutpoint $\hat z$ of age $a$ and origin $c \in \Crit$, we can
find $O_Z$  with boundary point $c$ such that $\pi(U_D) \cap
B_\eps(\pi(\hat z))$ intersects  $f^a(O_Z \cap B_\eps(c))$ for any
$\eps > 0$. Furthermore, there are rays $R_\phi$ and $R_{\phi'}$
(or possibly only one ray) landing at $c$ and intersecting
$\partial O_Z$ such that $f^a(R_\phi)$ and $f^a(R_\phi')$ land at
$\pi(\hat z)$ and intersect $\partial \pi(U_D)$.

There are distinct arrows leading to $D$ only if there is $D' \in
\D$ which is identified with $D$. But $D$ and $D'$ are only
identified if $\pi(D) = \pi(D')$, and the cutpoints, their origins
and ages coincide. Therefore the boundaries of $\pi(U_D)$ and
$\pi(U_{D'})$ are comprised of the same external rays together
with $\{ G(z) = 0\}$. It follows that $\pi(U_D) = \pi(U_{D'})$,
proving the lemma.
\end{proof}

\begin{lemma}
The system $(\hat U, \hat f)$ is Markov with respect to the
partition of canonical neighbourhoods, in the sense that if $\hat
f(U_D) \cap U_{D'} \neq \es$, then $\hat f(U_D) \supset
U_{D'}$.
\end{lemma}

\begin{proof}
This is a direct consequence of the previous proof.
\end{proof}

\section{Lifting measures}\label{sec_lifting}

In the previous section we introduced the Markov extensions and
canonical neighbourhoods for complex polynomials.  In this section
we will discuss the `liftability' properties of measures to the
Markov extension in the sense of Keller.  Our assumptions are
\eqref{SC} and \eqref{CO}.  We explain how they replace conditions
(2.2) and (2.3) of \cite{keller}.  This section gives the abstract
theory which is applicable to more general settings with this type
of Markov extension.  In subsequent sections the precise geometry
of $\J$, and thus the domains of $\hat\J$, play an important role
again.

Given a Borel $\sigma$--algebra $\B$ on $\J$ and a Borel
probability measure $\mu$ on $\J$, we will dynamically {\em lift}
this measure to a Borel probability measure $\hat \mu$ on $\hat
\J$. Our approach follows that of \cite{keller}.  We define a
method of obtaining $\hat \mu$ and then show that it is $\hat
f$--invariant and $\hat \mu \circ \pi^{-1} \ll \mu$.

We first introduce some notation.  For some space $X$, we let
$\cont_0(X)$ denote the set of continuous functions $\phi: X \to
\R$ with compact support.   For a set $A\subset X$, let $\chi_A: X
\to \{0,1\}$ be the characteristic function of $A$.

Let $i$ be the trivial bijection mapping $\J$ to $\hat \J_0$ (note
that $i^{-1}=\pi|_{\hat \J_0}$). Let
\begin{equation}
\hat\mu_0 \circ i = \mu\ {\rm and}  \ \hat{\mu}_n =
\frac{1}{n}\sum_{k=0}^{n-1} \hat\mu_0 \circ\hat{f}^{-k}.
\label{eqn:mulift}
\end{equation}

We will find some $\hat \mu$ to be a limit of a subsequence of
these measures.  Note that $\hat \J$ is in general not compact, so
the sequence $\{ \hat\mu_n \}_n$ may not have a subsequence with
limit in the weak topology.  Instead we use the vague topology,
see for example \cite{bil}. Given a topological space, we say that
a sequence of measures $\sigma_n$ converges to a measure $\sigma$
in the vague topology if for any function $\phi \in \cont_0(X)$,
we have $\lim_{n \to \infty} \sigma_n(\phi) = \sigma(\phi)$.  The
sequence $\{\hat\mu_n\}_n$ given in \eqref{eqn:mulift} has an
accumulation point in the vague topology.

\begin{definition}
A probability measure $\mu$ on $\J$ is {\em liftable} if a vague
limit $\hat \mu$ obtained in \eqref{eqn:mulift} is not identically
$0$.
\end{definition}

\begin{remark} Note that the measure $\mu \circ \pi$ on $\hat \J$ is
in general $\sigma$--finite, and not $\hat f$--invariant. The
lifted measure $\hat \mu$ distributes the mass of $\mu$ over the
domains of $\hat \J$ so as to become invariant, as we shall see
below.  Indeed $\mu$ is already called liftable if part of the
mass lifts to $\hat \J$; this is to accommodate non-ergodic
measures $\mu$. As in \cite{keller}, we generally wish to exclude
the possibility of $\hat \mu \equiv 0$ (where all the mass escapes
to infinity). Later we will find conditions to ensure that this
will not happen.
\end{remark}

The following theorem extends Theorem 2 of \cite{keller} from
ergodic invariant probability measures to general invariant
probability measures.

\begin{theorem}\label{thm:liftable}
Suppose that $\mu$ is an invariant probability measure on $\J$
satisfying \eqref{SC} and \eqref{CO}.  If $\hat\mu$ is a vague
limit point of any subsequence of measures given by
(\ref{eqn:mulift}) then it is $\hat f$--invariant and there is
some measurable function $0 \le \rho \le 1$ \st $\hat\mu \circ
\pi^{-1}= \rho\cdot\mu$. \label{thm:acip}
\end{theorem}

We will first state the theorem for ergodic invariant probability
measures, and then use the ergodic decomposition to generalise to
all invariant probability measures.

\begin{proposition}\label{prop:acip}
Suppose that $\mu$ is an ergodic invariant probability measure
satisfying \eqref{SC} and \eqref{CO}. If $\hat\mu$ is a vague
limit point of any subsequence of measures given by
(\ref{eqn:mulift}) and $\hat \mu \not\equiv 0$ then $\hat\mu$ is
an ergodic invariant measure and $\hat\mu \circ \pi^{-1}= \mu$.
\end{proposition}

Once we have shown that conditions (2.2) and (2.3) of
\cite{keller} can be replaced by \eqref{SC} and \eqref{CO} then
the proposition follows from \cite[Theorem 2]{keller}.

Theorem 1 of \cite{keller} implies that any ergodic invariant
probability measure $\mu$ can be lifted to a finite measure
$\hat\mu$ by applying \eqref{eqn:mulift}.  The conclusions of that
theorem also hold in our case. Conditions (2.2) and (2.3) of that
paper need to be assumed there in order to show that the lifting
process preserves ergodicity, and thus \cite[Theorem 2]{keller}
holds. The following lemma, which takes the role of \cite[Lemma
1]{keller}, shows that \eqref{SC} and \eqref{CO} are enough in our
case to draw the same conclusion here (i.e. lifting preserves
ergodicity and hence Proposition~\ref{prop:acip} holds).

Define $ \I := \left\{A \in \B: f^{-1}(A)=A\right\}  \mbox{ and }
\hat\I := \left\{A \in \hat\B: \hat f^{-1}(A)=A\right\} $

\begin{lemma} Let $\mu$ satisfy \eqref{SC} and \eqref{CO}.
Suppose that $\hat \mu$ is a vague limit of a subsequence of
$\{\mu_n\}_n$ \st $\hat \mu \circ\pi^{-1}=\mu$. Then $\pi^{-1}(\I)
= \hat \I \mod \hat \mu$. \label{lem:liftingergodicity}
\end{lemma}

\begin{proof} Suppose that $A \in \I$.
Then $\hat f^{-1} \circ \pi^{-1}(A) = \pi^{-1} \circ f^{-1}(A) =
\pi^{-1}(A)$ and so $\pi^{-1} (\I) \subset \hat \I$.

Conversely, suppose that $\hat A \in \hat\I$ and let $A = \pi(\hat
A)$.  Let $\hat B = \hat A \vartriangle \pi^{-1}(A)$.  We will
show that $\hat \mu(\hat B) = 0$. It follows from \eqref{SC} that
$\diam \ \hat Z_n[\hat z] \to 0$ for $\hat \mu$--a.e. $\hat z$.
Furthermore, \eqref{CO} implies that $\hat\mu_n(p) = 0$ for every
$n$ and $p\in\cut$. Therefore $\hat \mu(p) = 0$ as well.  Hence
$\hat\mu$--a.e. $z$ fulfilling the conditions of
Lemma~\ref{Lemfibersmeet} must be in case (a) of that lemma.
\footnote{ This is the same as saying that $(\hat\J,\hat
f,\hat\mu)$ satisfies condition (2.3) of \cite{keller}.
\label{rmk:SC CO}}

Therefore, for $\hat \mu$--a.e. $z_1 \in \hat B$, there exists
$z_2 \in \hat A$ and $n \ge 1$ \st $\hat f^n(z_1)= \hat f^n(z_2)$.
Hence $\hat\mu(\hat B)>0$ implies that $\hat A$ is not invariant;
a contradiction, whence $\hat \mu\left(\hat A \vartriangle
\pi^{-1}(A)\right) = 0$.  Thus, $\hat \I \subset \pi^{-1}(\I) \mod
\hat \mu$ and the lemma is proved.
\end{proof}

\begin{remark}
If $\hat\mu$ is an ergodic invariant probability measure on
$\hat\J$ such that $\hat\mu\circ\pi^{-1}=\mu$, then $\mu$ is
liftable. This is because it can be shown that for $\hat \mu_n$
defined as in \eqref{eqn:mulift}, $\hat\mu_n(\hat \J_R) \ge \hat
\mu(\hat J_R)$ for all $n, R \in \N$.   Moreover the lift of $\mu$
is absolutely continuous (and therefore equal) to $\hat \mu$.
Also, it follows from the proof of Theorem 2 in \cite{keller} that
given an ergodic invariant probability measure $\mu$ satisfying
\eqref{SC} and \eqref{CO}, there is at most one ergodic $\hat
f$--invariant probability measure $\hat\mu$ \st
$\hat\mu\circ\pi^{-1}=\mu$; so $\hat\mu$ is unique.

For liftable non--invariant measures, for example those considered
in Section~\ref{sec_conformal}, the measures
$\hat\mu\circ\pi^{-1}$ and $\mu$ are different.
\end{remark}

\begin{proof}[Proof of Theorem~\ref{thm:acip}]
Let $\B$ the $\sigma$-algebra of $\mu$-measurable sets, and let
\begin{equation}\label{eqn:ergdec}
\mu( \cdot ) = \int_Y \mu_y( \cdot ) \ d\nu(y)
\end{equation}
be the ergodic decomposition \label{page:erg_decomp} of $\mu$.
More precisely, the measure space $(Y,{\mathcal C},\nu)$ is used
to index the collection of all ergodic invariant probability
measures for $(\J, \B)$ and the probability measure $\nu$
satisfies \eqref{eqn:ergdec}. The diagram
\[
\begin{array}{ccc}
(\J,\B,\mu) & \stackrel{f}{\longrightarrow} & (\J,\B,\mu) \\[0.2cm]
\Pi \downarrow  & & \downarrow \Pi \\[0.2cm]
(Y, \mathcal C, \nu) &  \stackrel{\mbox{id}}{\longrightarrow} &
(Y, \mathcal C, \nu)
\end{array}
\]
commutes, the map $\Pi$ is such that $\Pi(z) = \Pi(z')$ if
$f^n(z)=f^m(z')$ for some $n,m \ge 0$, and $\mathcal C$ is the
finest $\sigma$-algebra such that $\Pi$ is $\B$-measurable. For
each $y \in Y$, $\Pi^{-1}(y)$ is called the {\em carrier} of
$\mu_y$; it is unique up to sets of $\mu_y$--measure $0$. For each
$y \in Y$, Proposition~\ref{prop:acip} states that there exists a
lifted measure $\hat \mu_y$ as the vague limit of $\{\hat
\mu_{y,n}\}_n$ constructed as in \eqref{eqn:mulift} (note that the
vague limit was independent of the subsequence chosen), and either
$\hat \mu_y \equiv 0$ or $\hat\mu_y(\hat \J) = 1$. Let $L = \{ y
\in Y : \mu_y \mbox{ is liftable}\}$.

{\bf Claim:} $L \in {\mathcal C}$; more precisely, there exists
$L' \in \B$ such that $\Pi:L' \to {\mathcal C}$ is well-defined
pointwise, and $\Pi(L') = L$.

To prove this claim, fix a countable $C^0$--dense subset $\hat\Phi
:= \{ \hat \phi_k \}_k$ of  $\cont_0(\hat\J)$ and a countable
collection of open intervals $\{ \mathcal U_l\}_l$ generating the
standard topology of $\R$. For each $y \in L$, we can use the set
$T_y$ of $\mu_y$--typical points as carrier. (Recall that $z$ is
called $\mu_y$--typical if the ergodic average $\frac1n
\sum_{i=0}^{n-1} \phi \circ f^i(z) \to \int \phi d\mu_y$ for each
continuous function $\phi:\J \to \R$.) As $\mu_y$ is liftable, the
lifted measure $\hat \mu_y$ has its set  $\hat T_y$ of
$\hat\mu_y$--typical points. Obviously $\hat T_y \subset
\pi^{-1}(T_y)$ and if $\hat z_0, \hat z_1 \in \pi^{-1}(z) \cap
\hat T_y$, then by Lemma~\ref{Lemfibersmeet}, there is $n$ such
that $\hat f^n(\hat z_0) = \hat f^n(\hat z_1)$. Therefore
$\hat\mu_y(\pi^{-1}(T_y) \vartriangle \hat T_y) = 0$, and a
fortiori, $i \circ \pi(\hat z) \in \hat T_y$ for each $\hat z  \in
\hat T_y$. Let $L'$ be the set of points $z \in \J$ such that
$i(z)$ is typical for $\hat \mu_y$ for some $y \in Y$. This is
exactly the set of points $z \in \J$ such that the ergodic
averages $\frac1n \sum_{j=0}^{n-1} \hat \phi_k \circ \hat f^j
(i(z))$  converge for each $k \in \N$, and at least one of the
limits $\neq 0$ (otherwise $z$ could only be typical for a
non--liftable measure $\mu_y \in {\mathcal C}$). Let
\[
X_{n,k,l} := \left\{ z \in \J : \frac1n \sum_{j=0}^{n-1}
\hat\phi_k \circ \hat f^j(i( z )) \in {\mathcal U}_l \right\},
\]
then
\[
L' = \left( \bigcap_{l \in \N} \bigcap_{k \in \N} \bigcup_{N \in
\N} \bigcap_{n \ge N} X_{n,k,l} \right) \cap \left( \bigcup_{ \{ l
\in \N \ : \ 0 \notin {\mathcal U}_l  \} } \bigcup_{k \in \N}
\bigcup_{N \in \N} \bigcap_{n \ge N} X_{n,k,l} \right).
\]
This set is obtained using countable operations on $\B$-measurable
sets $X_{n,k,j}$, so it belongs to the $\sigma$-algebra $\B$. This
proves the claim.

Let $\rho$ be the indicator function of $L'$.
Define
\[
\hat\mu(\hat A) := \int_Y \hat \mu_y \ d\nu(y) = \int_Y \int_A
\rho \circ \pi(\hat z) \ d\mu_y(\hat z) \ d\nu(y),
\]
whence $\hat \mu \circ \pi^{-1} = \rho \cdot \mu$. It remains to
show that $\hat \mu$ is the vague limit of the measures $\{\hat
\mu_n\}_n$ constructed in \eqref{eqn:mulift}.

Given $\eps > 0$, $\hat \phi\in \cont_0(\hat\J)$, for each $y \in
Y$ we can find $N = N(\eps,\hat\phi,y)$ such that
\[
|\hat\mu_{y,n}(\hat\phi) - \hat\mu_y(\hat \phi)|<\eps \mbox{ for
all } n \ge N.
\]
If $\mu_y$ is non--liftable, then $\hat \mu_y \equiv 0$; in this
case $|\hat \mu_{y,n}(\hat \phi)| < \eps$ for $n \ge N$.

Take $N_0$ so large that if $Y_0  = \{ y \in Y :
N(\eps,\hat\phi,y) > N_0\}$ then $\nu(Y_0) < \eps$. Then for $n
\ge N_0$,
\begin{eqnarray*}
|\hat \mu_n(\hat\phi) - \hat\mu(\hat\phi) | &\le& \int_{Y
\setminus Y_0} |\hat\mu_{y,n}(\hat\phi) - \hat\mu(\hat\phi) |
\ d\nu(y) \\
&& \hspace{2mm}
+ \int_{Y_0} |\hat\mu_{y,n}(\hat\phi) - \hat\mu(\hat\phi)| \ d\nu(y) \\
&\le& \int_{Y \setminus Y_0} \eps \ d\nu(y)
+ 2 \sup\hat\phi \ \nu(Y_0) \\
&\le & (1+2\sup\hat\phi) \eps.
\end{eqnarray*}
Since $\eps$ is arbitrary, $\hat \mu_n(\hat\phi) \to \hat
\mu(\hat\phi)$ as required.
\end{proof}

We will use the following lemma often in the forthcoming sections.

\begin{lemma}
Suppose that $\hat \mu$ is some measure on $\hat \J$ obtained from
applying (\ref{eqn:mulift}) to the probability measure $\mu$ on
$\J$.  If $\hat \mu\not\equiv 0$ and $\hat \mu\circ \pi^{-1} \ll
\mu$ then $\hat\nu:=\frac{\hat\mu}{\hat\mu(\hat\J)}$ is an
invariant probability measure on $\hat \J$. \label{lem:invmeas}
\end{lemma}

Note that the property $\hat \mu \circ \pi^{-1} \ll \mu$ is
immediate if $\mu$ is invariant. Indeed, in this case $\hat \mu_n
\circ \pi^{-1} = \mu$ for all $n$. So the lemma is useful when
$\mu$ is not invariant.

\begin{proof}
Let $\hat\phi\in \cont_0(\hat \J)$.  Define $\hat\nu_0:=\nu\circ
i^{-1}$. Then for any $R, n \ge 1$,
\begin{align*}
\hat\nu_n(\hat\phi\circ\hat f\cdot \chi_{\hat\J_R}) & =\frac 1 n
\sum_{j=0}^{n-1} \int_{\hat\J_R} \hat\phi\circ\hat f\
d(\hat\nu_0\circ \hat f^j) = \frac 1 n \sum_{j=0}^{n-1}
\int_{\hat\J_R} \hat\phi\circ\hat f^{j+1}\ d\hat\nu_0\\
&= \frac 1n\left( \sum_{j=0}^{n-1} \int_{\hat\J_R}
\hat\phi\circ\hat f^{j}\ d\hat\nu_0 + \int_{\hat\J_R}
\hat\phi\circ\hat f^{n}\ d\hat\nu_0 - \int_{\hat\J_R} \hat\phi\
d\hat\nu_0 \right)
\end{align*}
Therefore, $$\left|\hat\nu_n(\hat\phi\circ\hat f\cdot
\chi_{\hat\J_R})- \hat\nu_n(\hat\phi \cdot
\chi_{\hat\J_R})\right|< \frac{2\sup|\hat\phi|}{n}.$$ Letting $n,R
\to \infty$ we have proved the lemma.
\end{proof}

A dynamical system $(X, T, \mu)$ is said to be {\em dissipative}
if there is a {\em wandering set} of positive measure, i.e. a set
$A \subset X$ with $\mu(A) > 0$ \st $\mu( T^{-n}(A) \cap A)=0$ all
$n >0$. Otherwise the system is {\em conservative}. The system is
{\em totally dissipative} if there is no set $Y$ with $\mu (Y)>0$
and $\mu(T^{-1}(Y) \vartriangle Y) = 0$ \st $(Y, T,\mu|_Y)$ is
conservative.

In \cite{al}, the dissipativity/conservativity of various
quadratic polynomials with Feigenbaum combinatorics is
investigated.  For a lifted measure, we only see a conservative
part of the dynamics. This can be seen in the following lemma.

\begin{lemma}
Suppose that $(\J, f,\mu)$ is totally dissipative and $\hat\mu$ is
a measure obtained by applying \eqref{eqn:mulift} to the
probability measure $\mu$. If $\hat\mu\circ\pi^{-1} \ll \mu$ then
$\hat\mu \equiv 0$. \label{lem:totdiss}
\end{lemma}

\begin{proof}
We start with the following claim: no measure $\hat \mu$ obtained
by applying (\ref{eqn:mulift}) to a probability measure $\mu$ can
have wandering sets of positive $\hat\mu$--measure. Indeed,
suppose that $\hat A \subset \hat\J$ has $\hat f^{-n} (\hat A)
\cap \hat A = \es$ for all $n \ge 1$. Suppose that $\hat \mu$ is a
vague limit of $\{\hat \mu_{n_k}\}_k$. We will show that $\hat \mu
(\hat A)=0$. Note that $\sum_{i=0}^{n-1} \hat\mu_0 \left(\hat
f^{-i}(\hat A)\right) \le 1$ since the domains $\hat f^{-i}(\hat
A)$ are disjoint.  Thus
\[
\hat \mu (\hat A) = \lim_{k \to \infty} \frac{1}{n_k}
\sum_{i=0}^{n_k-1} \hat\mu_0 \left(\hat f^{-i}(\hat A)\right)=0.
\]
This proves the claim, and hence $\hat \mu$ is conservative.

Now suppose that $\hat \mu \not\equiv 0$, then $\hat \mu \circ
\pi^{-1}$ is an $f$--invariant measure which can be normalised,
say $\mu_0 := \frac1{\hat \mu(\hat \J)} \ \hat \mu \circ
\pi^{-1}$. Let $Y$ be the carrier of $\mu_0$ (or more precisely,
the union of the carriers of all liftable ergodic measures present
in the ergodic decomposition of $\mu_0$), then since $\hat \mu
\circ \pi^{-1} \ll \mu$, $\mu(Y) > 0$ and $\mu(f^{-1}(Y)
\vartriangle Y) = 0$. Since $\mu$ is totally dissipative, there
must be a wandering set $A \subset Y$, and hence $\pi^{-1}(A)$ is
wandering for $\hat \mu$. This contradiction proves the lemma.
\end{proof}

\section{Inducing}\label{sec_induce}

A particularly useful property of Markov extensions is that they
easily enable one to construct uniformly expanding induced systems
with bounded distortion, provided the measure $\mu$ is liftable.
In fact, any first return map on the Markov extension corresponds
to an induced (jump) transformation of the original system, and
under mild conditions, the reverse is true as well, cf.
\cite{bruin1}. If $\hat W \subset \hat \J$, let us write $\hat
F_{\hat W}$ for the first return map to $\hat W$, i.e. $\hat F(z)
= \hat f^{\tau(z)}(z)$ where $\tau = \tau_{\hat W}:{\hat W} \to
\N$ is the first return time to $\hat W$. Let $\tau^n(z)$ denote
the $n$--th return time, i.e. $\tau^1(z) = \tau(z)$ and $\tau^n(z)
= \tau^{n-1}(z) + \tau(\hat F^{n-1}(z))$. For our purposes, we are
most interested in subsets $\hat W$ of some domain $D\in \D$ that
are bounded away from the cutpoints of $D$. As a result, any such
set $\hat W$ is compactly contained in the canonical neighbourhood
$U_D$ of $D$, and by the Markov property of $(\sqcup_D U_D, \hat
f)$, any branch of $\hat F^n_{\hat W}=\hat f^{\tau^n({\hat
W}_0)}:{\hat W}_0 \to {\hat W}$ for any $n \in \N$ is extendible
to a univalent onto map $\hat f^{\tau^n}:V_0 \to U_D$.

Given $\delta  > 0$ and $M > 0$ we say that $z$ {\em reaches large
scale at time $j$} if there are neighbourhoods $\C \supset V_0
\supset V_1 \owns z$ such that $f^j:V_0 \to f^j(V_0)$ is
univalent, $f^j(V_1)$ contains a round ball of radius $\delta$
(measured in Euclidean distance) and $\modulus(V_0,V_1)
> M$, see \cite{Mil} for definitions. It follows from the Koebe
distortion theorem, see \cite[Theorem 1.3]{Pom}, that there exists
$K = K(M)$ such that the distortion
\[
\dist(f^j|_{V_1}) := \sup_{z,z' \in V_1} \frac{|Df^j(z)|}{|Df^j(z')|}
\le K.
\]

\begin{lemma}\label{lem:large_scale}
Let $\mu$ be an ergodic $f$--invariant probability measure
satisfying \eqref{SC} and \eqref{CO}. Then $\mu$ is liftable if
and only if there exist $\delta
> 0, v > 0$ and $M > 0$ such that for $\mu$--a.e. $z \in \J$
\begin{equation}
\liminf_{n \to \infty} \frac1n \#\{ 0 \le j < n : z \mbox{ reaches
large scale for } \delta, M \hbox{ at time } j\} \ge v.
\label{eq:large scale}\end{equation} In this case $\hat \mu \circ
\pi^{-1} =\mu$.
\end{lemma}

\begin{proof}
First assume that $\mu$ is liftable and let $\hat \mu$ be the
lifted measure.  Let $D\in\D$ and let $\hat W \subset U_D$ be an
open set bounded away from the cutpoints of $D$ such that (using
\eqref{CO}) $v := \hat \mu({\hat W}) > 0$. Take $\delta$ such that
$U_D$ and $\hat W$ contain round balls of radius $\delta$. Also
$\hat W$ is compactly contained in $U_D$, so $M :=
\modulus(U_D,{\hat W}) > 0$.  Let $z$ be a typical point for $\mu$
and let $\hat z =i(z)$. By Birkhoff's Ergodic Theorem,
\[
\lim_{n \to \infty} \frac1n \#\{ 0 \le j < n : \hat f^j(\hat z)
\in {\hat W} \} = \hat\mu({\hat W}).
\]
By the Markov property $z$ reaches large scale for $\delta, M$ at
time $j$ if $\hat f^j(\hat z) \in {\hat W}$. It follows that
$\pi(\hat z)$ has reached large scale at time $i$ as well and so
the first implication follows.

Conversely, suppose that $\mu$--a.e. $z$ satisfies \eqref{eq:large
scale}.  We say that $z\in H_R$ if,  given $\hat z$ \st $\pi(\hat
z)=z$ and $B_\delta(\hat z)$ has $\modulus(U_D, B_\delta) > M$,
then $\pi^{-1}(Z_R[z])$ contains no cutpoint of $D$. By assumptions
\eqref{SC} and \eqref{CO}, $\mu(H_R) \to 1$ as $R \to \infty$, say
$\mu(H_R) > 1-\eta(R)$ where $\lim_{R \to \infty} \eta(R) = 0$.

If $z$ reaches large scale for $\delta > 0$ and $M>0$, at iterate
$j$, then for $\hat z  = i(z)$ and the domain $D \owns \hat
f^j(\hat z)$, $\hat Z_R[\hat f^j(\hat z)]$ contains no cutpoint of
$D$. It follows that $\hat f^{j+R}(\hat z) \in \hat \J_R$.
Therefore, given $\eps > 0$ and a $\mu$--typical point $z$, there
exists $n_0(z)$ such that for $n \ge n_0(z)$
\begin{align*}
(i) & \hspace{10mm} \frac1n \#\{ 0 \le j <
n : f^j(z) \in H_R \} > 1-2\eta(R), \\[0.2cm]
(ii) & \hspace{10mm} \frac1n \#\{ 0 \le j < n : \hat f^j(\hat z)
\in \hat \J_R\}
> \frac{v}{1+\eps}.
\end{align*}
Take $M$ so large that $n_0(z) \le M$ for all $z$ in a set of
$\mu$--measure $\ge 1-\eps$. Then
\[
\hat \mu_n(\hat \J_R) = \frac1n \sum_{j=1}^{n-1} \hat \mu_0\circ
\hat f^{-j}(\hat \J_R) \ge v \left(\frac{r-1}{r}\right)
\left(\frac{1-\eps}{1+\eps}\right) (1-2 \eta(R)),
\]
for all $n \ge rM$. As $r \in \N$ and $\eps > 0$ are arbitrary,
any vague limit point of $\{ \hat \mu_n\}_n$ satisfies $\hat
\mu(\hat \J_R) \ge v (1-2\eta(R))$, which is positive for $R$
sufficiently large. By Proposition~\ref{prop:acip} this means that
the ergodic measure $\mu$ is liftable and $\hat\mu\circ\pi^{-1}
=\mu$.
\end{proof}

\begin{remark}
Our notion of `reaching large scale' with positive frequency is
stronger than the notion of {\em induced hyperbolicity} in \cite
{GS1}. Note also that in fact the proof above shows that given
$\delta, M>0$, if\textsl{}
\[
\liminf_{n \to \infty} \frac1n \#\{ 0 \le j < n : z \mbox{ reaches
large scale for } \delta, M \hbox{ at time } j\}
> 0.
\]
on a positive measure set for {\em any} probability measure $\mu$
then the measure $\hat \mu$ obtained from (\ref{eqn:mulift}) is
non--zero. \label{rmk:pos freq}
\end{remark}

If a point $\hat z$ visits a compact part $\hat \J_R$ of the
Markov extension with positive frequency, then the majority of
these visits are at a certain distance away from cutpoints in
$\hat \J_R$. This is made precise in the following lemma. As a
consequence, $z = \pi(\hat z)$ will go to large scale (with
bounded distortion) with positive frequency.

\begin{lemma} \label{lem:visits_to_cutpoints}
Suppose that \eqref{CO} is satisfied.  For each domain $D\in\D$
and $\eps> 0$, there exists $\delta> 0$ such that if
$\hat X = \bigcup_{p\in \cut\cap D} B_\delta(p)$, then for every
invariant probability measure $\hat \mu$ on $\hat \J$, $\hat
\mu(\hat X) < \eps$.
\end{lemma}

\begin{proof}
Suppose that $D\in \D$ has $\lev(D)=n$ and $p\in \cut_D := \cut
\cap D$. Then $p$ has age $m\le n$.  Since the domain containing
$\hat f^{n-m+j}(p)$ must have level at least $j+n$, $p$ can return
to $D$ under iteration by $\hat f$ a maximum of $n-m$ times.
Therefore, there exists $n_0 \ge 1$ \st $f^{n_0+k}(\cut_D) \notin
D$ for all $k \ge 1$. So for any $j_0$, there exists $\delta>0$
\st $\hat f^{j+n_0}(B_\delta(p)) \cap D = \es$ for $1 \le j
\le j_0$. Take  $\hat X := \bigcup_{p\in \cut_D} B_\delta(p)$. If
$k \in \{0,\dots, j_0\}$ and $l > k$ is such that $\hat
f^{-k}(\hat X) \cap f^{-l}(\hat X) \neq \es$, then $\hat
f^{l-k}(\hat X) \cap D \neq \es$. Therefore there are at
most $2n_0$ numbers $l \in \{ 0, \dots, n_0+j_0\}$ such that $\hat
f^{-k}(\hat X) \cap f^{-l}(\hat X) \neq \es$. Furthermore
$\hat\mu(\hat f^{-k}(\hat X)) = \hat \mu(\hat X)$. It follows by
Lemma~\ref{lem:invmeas} that $1 \ge \hat \mu( \cup_{k=0}^{j_0}
\hat f^{-k}(\hat X) ) \ge \frac{j_0}{2n_0+1} \hat\mu(\hat X)$. To
complete the proof, take $j_0 > 2n_0/\eps$ and get $\delta>0$
accordingly.
\end{proof}

\section{Liftability and positive Lyapunov exponents}
\label{sec_lyap}

Given a dynamical system $(X, g, \nu)$, let $\phi_g = \log |Dg|$
wherever this is defined and $\lambda_g(\nu) = \int \phi_g \ d\nu$
be the Lyapunov exponent of $\nu$. The pointwise (upper and lower)
Lyapunov exponents at a point $x \in X$ are denoted as
$\lambda_g(x)$ (and $\overline \lambda_g(x)$ and $\underline
\lambda_g(x)$ respectively) wherever these are well-defined.

\begin{proposition} \label{prop:lyap equal}
Suppose that \eqref{CO} holds.  If $\mu$ is an ergodic invariant
liftable probability measure, with lifted measure $\hat \mu$, then
\begin{enumerate}
\item[(a)] $\phi_f$ is integrable with respect to $\mu$;
\item[(b)] $\lambda_f(\mu) = \lambda_{\hat f}(\hat \mu) > 0$.
\end{enumerate}
\end{proposition}

\begin{proof}[Proof of Proposition~\ref{prop:lyap equal}]
Note that $\hat \mu$ is an invariant measure: for example see
Lemma~\ref{lem:invmeas}.  By Lemma~\ref{lem:visits_to_cutpoints},
we may take domain $D\in\D$ and $\hat W \in D\cap\hat\P_n$ such
that $\hat W$ is compactly contained in $U_D$ and $\hat \mu({\hat
W}) > 0$. By the Poincar\'e Recurrence Theorem, $\hat F_{\hat
W}:\bigcup_j\hat W_j \to \hat W$, the first return map to $\hat
W$, is defined $\hat \mu$--a.e.  Given $z\in \hat W_j\sm \bd \hat
W_j$ there is an open neighbourhood $U$ of $x$ \st $\hat F_{\hat
W}$ extends to this neighbourhood.  In particular $D\hat F_{\hat
W}$ is defined for all $z\in \bigcup_j\hat W_j\sm \bd \hat W_j$.
In particular, since $\hat\mu(\bd \hat W_j)=0$ (otherwise
\eqref{CO} is contradicted), the derivative is defined for
$\hat\mu$ a.e. $z\in \hat W$.  Each branch of $\hat F^n_{\hat W}$
is extendible to $U_D$, so by the Koebe distortion theorem,
$\kappa := \inf\{ |D\hat F(z)| : \hat F(z) \mbox{ is
well--defined} \}  > 0$. In fact, there is $N$ such that $\inf\{
|D\hat F^N(z)| : \hat F^N(z) \mbox{ is well--defined} \} \ge 2$
(one consequence of this is given in Remark~\ref{rmk:parab}
below).

The measure $\hat \mu_{\hat W} := \frac{1}{\hat\mu({\hat W})}\hat
\mu|_{\hat W}$ is an $\hat F_{\hat W}$--invariant probability
measure, and Kac's Lemma implies that
\[
\int \tau \ d\hat\mu_{\hat W} =  \frac{1}{\hat\mu({\hat W})} <
\infty,
\]
where $\tau = \tau_{\hat W}$ is the first return time by $\hat f$
to $\hat W$. Moreover $D\hat F^n(\hat z) = D\hat f^{\tau^n(z)}(z)$
and if $\hat z$ is typical for $\tau$, then denoting $L_f=\sup_{z
\in \J} |Df(z)|$,
\begin{eqnarray*}
0 &<& \hat \mu({\hat W}) \liminf_{n \to \infty}
\frac1n \log |D\hat F^n(\hat z)| \\
&=& \hat \mu({\hat W}) \lim_{n \to \infty} \frac{\tau^n(\hat z)}n
\liminf_{n \to \infty}\frac1{\tau^n(\hat z)}
\log |D\hat f^{\tau^n(\hat z)}(\hat z)| \\
&\le& \overline \lambda_{\hat f}(\hat z) \le \log L_f < \infty.
\end{eqnarray*}
For $L<\infty$, take $\hat\Phi_L = \min\{ L, \log|D\hat F(z)|\}$.
Then $\hat\Phi_L$ is bounded and hence $\hat\mu_{\hat
W}$--integrable, and for $\hat\mu_{\hat W}$--a.e. $\hat z$
\begin{eqnarray*}
0 &<& \int_{\hat W} \hat\Phi_L\ d\hat\mu_{\hat W}
= \lim_{n \to \infty} \sum_{k=0}^{n-1} \hat\Phi_L(\hat F^k(\hat z))\\
&\le& \lim_{n \to \infty} \sum_{k=0}^{n-1} \log |D\hat F(\hat
F^k(\hat z))|
= \lim_{n \to \infty} \frac1n \log |D\hat F^n(\hat z)| \\
&\le& 
\frac{1}{\hat\mu({\hat W})} \log L_f < \infty.
\end{eqnarray*}
The Monotone Convergence Theorem gives that $\log |D\hat F| =
\lim_{L \to \infty} \hat\Phi_L$ is $\hat\mu_{\hat W}$--integrable
and
\[
\int_{\hat W} \log |D\hat F| \ d\hat\mu_{\hat W}  =  \lim_{n \to
\infty} \frac 1n\log |D\hat F^n(z)| \qquad \mbox{ $\hat\mu_{\hat
W}$--a.e.}
\]
We can apply the same argument to $\phi_U := \max\{ -U, \log |Df|
\}$, which is $\mu$--integrable: for $\mu$--a.e. $z \in W :=
\pi(\hat W)$ and $\hat z \in \hat W$ such that $\pi(\hat z) = z$
\begin{eqnarray*}
\log L_f &\ge& \int \phi_U \ d\mu = \lim_{n \to \infty}
\frac1{\tau^n(\hat z)}
\sum_{k=0}^{\tau^n(\hat z)-1} \phi_U(f^k(z)) \\
&=& \lim_{n \to \infty} \frac{n}{\tau^n(\hat z)} \lim_{n \to
\infty} \frac{1}{n} \sum_{k=0}^{n-1} \ \ \sum_{j = \tau^k(\hat
z)}^{\tau^{k+1}(\hat z)-1}
\phi_U(f^j(z)) \\
&\ge& \hat\mu({\hat W}) \lim_{n \to \infty} \frac1n
\sum_{k=0}^{n-1}
\log |D\hat F(\hat F^k(\hat z))| \\
&=& \hat\mu({\hat W}) \int_{\hat W} \log |D\hat F| \
d\hat\mu_{\hat W}
> 0.
\end{eqnarray*}
The Monotone Convergence Theorem implies that $\log |Df| = \lim_{U
\to \infty} \phi_U$ is $\mu$--integrable. Hence $\underline
\lambda_f(z) = \overline \lambda_f(z) = \lambda_f(z)$ for
$\hat\mu$--a.e. $z$ and
\[
\lambda_f(z) = \lambda_{\hat f}(z) = \hat\mu({\hat W})
\lambda_{\hat F}(z) \ge \hat\mu({\hat W}) \frac{\log 2}{N}.
\]
\end{proof}

\begin{remark} Our set-up of dendrite Julia sets necessarily excludes
the existence of neutral periodic cycles, but also when the
construction is extended to more general Julia sets, for example
with Siegel disks or Leau--Fatou petals (cf. \cite{Mil}), the
proof of this proposition shows that Dirac measures on parabolic
periodic points are not liftable to the Markov extension.
\label{rmk:parab}
\end{remark}

In the next result, let $\hat W$ and $\hat F$ be as in the proof
of Proposition~\ref{prop:lyap equal}.

\begin{proposition}
If \eqref{SC} and \eqref{CO} hold and $\mu$ is invariant, ergodic
and liftable, then
\[
h_\mu(f) = h_{\hat\mu}(\hat f) = \hat\mu(\hat W) h_{\hat\mu_{\hat
W}}(\hat F).
\]
\end{proposition}

\begin{proof}
The first equality can be shown in the same way as  Theorem 3 from
\cite{keller}. Note that \eqref{SC} by itself does not imply that
the partition $\P_1$ generates the Borel $\sigma$-algebra; the
condition used by Keller. But Keller's proof relies on the
Shannon--McMillan--Breiman Theorem, which only uses that $Z_n[z]
\to 0$ $\mu$-a.e., which is indeed condition \eqref{SC}. Otherwise
this equality follows from the fact that a countable--to--one
factor map preserves entropy, provided the Borel sets are
preserved by lifting, see \cite{DS}.

The second inequality is Abramov's formula, see \cite{Ab}.
\end{proof}

Given an invariant probability measure $\mu$ on $(\J, f)$, let
$(\tilde \J, \tilde f, \tilde \mu)$ be the natural extension. Each
$\tilde z \in \tilde\J$ is represented as a sequence $(\tilde z_0,
\tilde z_1, \dots)$ such that $\tilde z_j = f^j(\tilde z)\in \J$,
and $\tilde f(\tilde z_0, \tilde z_1, \dots) = (f(\tilde z_0),
\tilde z_0, \tilde z_1 \dots)$. Define $\tilde \pi_k(\tilde z) :=
\tilde z_k\in \J$.

The motivation for the following lemma, and the idea for the next
theorem are based on a result of \cite{ledinterval}.  For the
remainder of this section, we will always consider subsets of
$\J$, and we will use the relative topology on $\J$.

\begin{lemma}\label{lem:to_large_scale}
Suppose that $\mu$ is an invariant probability measure satisfying
\eqref{SC}. Given $\tilde z \in \tilde \J$, let $W(\tilde z) = \cap_{k
\ge 1} f^k Z_k[\tilde z_k]$ and $r(\tilde z) = \sup\{ \rho \ge 0;
B_{\rho}(\tilde z_0) \subset W(\tilde z)\}$. If $r(\tilde z)
> 0$ for $\tilde\mu$--a.e. $\tilde z$, then applying
\eqref{eqn:mulift} to $\mu$ gives rise to an invariant probability
measure $\hat \mu$ on $\hat \J$.
\end{lemma}

\begin{proof}
Let $\eps > 0$ be arbitrary and $r_0 > 0$ be such that $\tilde
\mu(A) > 1-\eps$ for $A = \{ \tilde z \in \tilde \J : r(\tilde z)
> r_0\}$. Then for each $k \ge 0$, $A_k := \tilde \pi_k(A)$
satisfies $1-\eps \le \mu(A_k) \le \mu \circ f^{-k}(A_0)$, and for
each $z \in A_k$, $B_{r_0}(f^k(z)) \subset f^k(Z_k[z])$.

Given $\hat z\in \hat\J$, define $D_{\hat z}\in \D$ to be the
domain containing $\hat z$. Take $R$ so large that if $\hat z \in
\hat \J$ and $B_{r_0}(\hat z)$ contains no cutpoints of $D_{\hat
z}$ then $\hat Z_R[\hat z] \subset B_{r_0}(\hat z)$. Define $K_R
:= \{ \hat z \in \hat \J : \hat Z_R[\hat z] \mbox{ contains no
cutpoint of }D_{\hat z}\}$.  Note that if $\hat z \in K_R$ then
$\hat f^R(\hat z)\in \hat\J_R$. If $z \in A_k$, then $f^k(Z_k[z])
\supset B_{r_0}(f^k(z))$ and by the Markov property of $(\hat \J,
\hat f)$, letting $\hat z:=i(z)$ we have $\hat f^k(\hat Z_k[\hat
z]) = D_{\hat f^k(\hat z)} \supset B_{r_0}(\hat f^k(\hat z))
\supset \hat Z_R[\hat f^k(\hat z)]$.
It follows that $\hat f^k(\hat z) \in K_R$, and therefore $\hat
f^{k+R}(\hat z) \in \hat \J_R$. This shows that
\[
\hat\mu_0 \circ \hat f^{-(k+R)}(\hat \J_R) \ge \hat\mu_0 \circ
\hat f^{-k}(K_R) \ge \hat\mu_0\circ i(A_k) = \mu(A_k) \ge 1-\eps
\]
for all $k \ge 0$. Therefore any vague limit point $\hat\mu$ of
$\{ \hat \mu_n\}_n$ satisfies $\hat\mu(\hat\J_R) \ge 1-\eps$, and
because $\eps > 0$ was arbitrary, $\hat\mu \not\equiv 0$ and in
fact $\hat\mu(\hat \J) = 1$.
\end{proof}

\begin{theorem}
Let $\mu$ be an $f$--invariant probability measure
satisfying \eqref{SC} and \eqref{CO}, such that $\lambda_f(z)>0$ $\mu$--a.e.
Then $\mu$ is liftable and $\hat\mu\circ\pi^{-1}\ll \mu$.
\label{thm:liftable from pos lyap}
\end{theorem}

\begin{proof}
We can apply \cite[Theorem 3.17]{EL}, which says that in this
setting there exists a partition $\eta$ of $\tilde\J$ into open
sets (recall we are using the relative topology on $\tilde\J$
here) \st for $\tilde\mu$--a.e. $\tilde z$, $\tilde f$ has bounded
distortion on the element $\eta(\tilde z)$ of $\eta$
containing $\tilde z$. More precisely,
there is constant $K(\tilde z) \ge 1$ \st
\[
K(\tilde z)^{-1} \le
\frac{|Df^n(\tilde x_n)|}{|Df^n(\tilde y_n)|} \le K(\tilde z)
\]
for all $n \ge 0$ and components $\tilde x_n, \tilde y_n$ of
$\tilde x, \tilde y \in \eta(\tilde x)$.
Define $r(\tilde z) := \sup\{
\rho \ge 0 : B_{\rho}(\tilde z_0) \subset \pi(\eta(\tilde z))\}$.
This is a $\tilde\mu$--measurable function which is strictly
positive $\tilde\mu$--a.e.

The general idea behind this result is from Ruelle, see
\cite{Rue}, and is usually presented as a `local unstable manifold
theorem'. It is given in the complex setting in \cite{ledrational}
and is discussed in \cite{EL}.  An alternative proof is presented
in Section 9 of \cite{PU}.

Now notice that for any $\tilde z \in \tilde\J$, $\eta(\tilde
z)\subset W(\tilde z)$ (as defined in Lemma~\ref{lem:to_large_scale}),
otherwise $\eta(\tilde z)\cap\bd
f^k Z_k[\tilde z_k] \neq \es$ for some $k$, 
which implies that distortion is
unbounded; a contradiction. Let $r'(\tilde z) := \sup\{ \rho \ge 0
: B_{\rho}(\tilde z_0) \subset W(\tilde z) \}$, then $r'(\tilde z)
\ge r(\tilde z)> 0$, $\tilde \mu$--a.e. By
Lemma~\ref{lem:to_large_scale}, applying (\ref{eqn:mulift}) to
$\mu$ gives a measure $\hat \mu \not\equiv 0$.
It follows from Theorem~\ref{thm:liftable} that $\hat \mu \circ
\pi^{-1} \ll \mu$ and $\hat\mu$ is invariant.
\end{proof}

\begin{corollary}
Let $\mu$ be an invariant probability measure satisfying
\eqref{SC} and \eqref{CO}. If the measure theoretic entropy
$h_{\mu}(f) > 0$, then $\mu$ is liftable.
\end{corollary}

\begin{proof}
It follows from the Ruelle inequality \cite{PU} that for ergodic invariant
measures, $\nu$, $h_{\nu}(f) \le 2\lambda(\nu)$.  Therefore, if we consider the ergodic decomposition,
our assumption implies that there is a positive $\mu$--measure set of $z$ with
$\lambda_f(z)>0$. Thus Theorem~\ref{thm:liftable from pos lyap}
implies that $\mu$ is liftable. (In the interval case, Keller
\cite{keller} gave a proof based on a counting argument of paths
high up in the tower. This type of proof can be used here too; see
the appendix for our counting argument, which is to be used in the
next section.)
\end{proof}

\section{Conformal measure}\label{sec_conformal}

In this section we discuss the liftability properties of conformal
measure.
Sullivan \cite{sull} showed that all rational maps on the Riemann
sphere have a conformal measure for at least one minimal $\delta
\in (0,2]$. We would like to emphasise that $\mu_\delta$ is not
invariant, but when $\mu_\delta$ is liftable, then $\hat
\mu_\delta$ (normalised) projects to an invariant probability
measure, say $\nu = \alpha \cdot \hat \mu_\delta \circ \pi^{-1}$,
where $\alpha \ge 1$ is the normalising constant.

Our first lemma is that $\nu$ is absolutely continuous,
generalising Proposition~\ref{prop:acip} to $\delta$-conformal
measure. It can be expected that this lemma generalises to other
non--invariant probability measures too, provided there is
distortion control.

\begin{lemma}\label{lem:conf_ac}
Suppose that a conformal measure $\mu_\delta$ on $\J$ satisfies
\eqref{CO}. Let $\hat \mu$ be a measure on $\hat \J$ obtained as a
vague limit of \eqref{eqn:mulift}. Then $\hat \mu \circ \pi^{-1}
\ll \mu_\delta$.
\end{lemma}

\begin{proof}
We suppose that $\hat \mu_\delta \not\equiv 0$, otherwise there is
nothing to prove.  Suppose that $\hat\mu_{n_k} \to \hat\mu_\delta$
as $k \to \infty$.

If the lemma is not satisfied then there exists $\eps>0$ and a set
$\hat A \subset \hat \J$ which has $\hat \mu_\delta (\hat A)>
\eps$, but $\mu_\delta (A) =0$ for $A=\pi(\hat A)$. We may assume
that $\hat A$ is contained in some domain $D \in\D$.

Due to \eqref{CO}, we can assume that $\hat A$ is compactly
contained inside $U_D$. Therefore $\hat
\mu_\delta\left(U_D\right)>0$. Choose some $\hat B$ compactly
contained in $U_D$ with $\mu_\delta(B)>0$ where $B=\pi(\hat B)$.
We take some neighbourhood $U$ containing both $A$ and $B$ which
is compactly contained in $U_D$. There is some $C>0$ \st for any
$x, y \in U$, for each branch of the inverse map we have
\[
\left|\frac{D\hat{f}^{-n}(x)}{D\hat{f}^{-n}(y)}\right|<C \hbox{
for all } n \ge 1.
\]
Supposing that $\delta>0$ is the exponent of the conformal
measure, we have
\[
\frac{\hat\mu_\delta(\hat A)}{\hat\mu_\delta(\hat B)} =
\frac{\lim_{k \to \infty} \frac{1}{n_k} \sum_{j=0}^{n_k-1}
\hat\mu_0 (\hat f^{-j}(\hat A))}{\lim_{k \to \infty} \frac{1}{n_k}
\sum_{j=0}^{n_k-1} \hat\mu_0 (\hat f^{-j}(\hat B))} \le
{C^{2\delta}}\ \frac{\mu_\delta(A)}{\mu_\delta(B)}.
\]
But while the left hand side is positive, the right hand side is
$0$, so we have a contradiction.  Thus we obtain absolute
continuity as required.
\end{proof}

Combining Remark~\ref{rmk:pos freq} and Lemmas~\ref{lem:invmeas}
and \ref{lem:conf_ac}, we get the following corollary.

\begin{corollary}
Suppose that $\mu$ is conformal measure $\mu_\delta$ satisfies
\eqref{SC} and \eqref{CO}. If for given $\delta, M>0$,
\[
\liminf_{n \to \infty} \frac1n \#\{ 0 \le j < n : z \mbox{ reaches
large scale for } \delta, M \hbox{ at time } j\}> 0,
\]
then $\mu_\delta$ is liftable.
\end{corollary}


The following lemma and theorem are similar to part of the
statement of Theorem B in \cite{ledrational}.  We supply a proof
for completeness.

\begin{lemma}\label{lem:conformal equiv}
Assume that conformal measure $\mu_\delta$ satisfies \eqref{SC}
and \eqref{CO}. If $\mu_\delta$ is liftable, and $\nu = \hat
\mu_\delta \circ \pi^{-1}$, then $\mu_\delta$ and $\nu$ are
equivalent.  Moreover, $\mu_\delta$ and $\nu$ are ergodic.
\end{lemma}

\begin{proof}
It was shown in Lemma~\ref{lem:conf_ac} that $\nu \ll \mu_\delta$.
Let us prove that $\psi:= \frac{d\nu}{d\mu_\delta}$ is a positive
density.

Let $A \subset J$ be an open set such that $\nu(A) > 0$.  Let
$\hat\nu$ be the measure obtained from applying \eqref{eqn:mulift}
to $\nu$. As $\nu = \hat \nu \circ \pi^{-1}$, there must be some
$D\in\D$ such that $\hat \nu(D \cap \pi^{-1}(A)) > 0$. By
replacing $A$ by an appropriate cylinder set $Z \in \P_n$, we can
assume (using \eqref{SC} and \eqref{CO}) that $\hat Z :=
\pi^{-1}(Z) \cap D$ is bounded away from the cutpoints of $D$ and
such that $\hat \nu(\hat Z) > 0$. Moreover, as $\hat Z \in \hat
\P_n$, no boundary point of $\hat Z$ (relative to $D$) returns to
$\hat Z$.

Let $\hat F$ be the first return map to $\hat Z$.  We will show that
we may apply
the Folklore Theorem to this map.  First note that if $\hat Z$ is
chosen sufficiently small, then all the branches of $\hat F$ are
expanding and the Koebe Lemma implies that they have bounded
distortion.  Let $\hat B\subset \hat Z$ be the set of points in
$\hat Z$ which never return to $\hat Z$. We now wish to check that
$\mu_\delta\circ\pi(\hat B)=0$.  We use the same technique as in
the proof of Lemma~\ref{lem:conf_ac}. By Poincar\'e recurrence we
know that $\hat\mu_\delta(\hat B)=0$. We let $B:=\pi(\hat B)$ and
suppose that $\mu_\delta(B)>0$, and will show this leads to a
contradiction. As in the proof of Lemma~\ref{lem:conf_ac}, we can
use a distortion argument to show that
\[
0=\frac{\hat\mu_\delta(\hat B)}{\hat\mu_\delta(\hat Z)} =
\frac{\lim_{k \to \infty} \frac{1}{n_k} \sum_{j=0}^{n_k-1}
\hat\mu_0 (\hat f^{-j}(\hat B))}{\lim_{k \to \infty} \frac{1}{n_k}
\sum_{j=0}^{n_k-1} \hat\mu_0 (\hat f^{-j}(\hat Z))} \ge \frac
1{C^{2\delta}}\ \frac{\mu_\delta(B)}{\mu_\delta(Z)}.
\]
But since the right hand side is bounded away from zero, we have a
contradiction.

We can now apply the Folklore Theorem to
$\frac{1}{\mu_\delta(\pi(\hat Z))} \mu_\delta \circ \pi|_{\hat
Z}$, which yields an ergodic $\hat F$--invariant probability
measure $\hat \nu_{\hat Z}$, with density $\hat \psi =
\frac{d\hat\nu_{\hat Z}}{d\mu_\delta \circ \pi|_{\hat Z}}$ bounded
above and bounded away from zero. Since $\mu_\delta$ is liftable
and, by Lemma~\ref{lem:conf_ac}, the lifted measure $\hat \nu$
satisfies $\hat \nu \circ \pi^{-1} \ll \mu_\delta$ we have
$\frac1{\hat \nu(\hat Z)} \hat \nu|_{\hat Z}\ll \hat\nu_{\hat Z}$.
Since $\hat\nu_{\hat Z}$ is ergodic and both $\hat\nu_{\hat Z}$
and $\frac1{\hat \nu(\hat Z)} \hat \nu|_{\hat Z}$ are $\hat
F$--invariant probability measures, $\hat\nu_{\hat Z} =
\frac1{\hat \nu(\hat Z)} \hat \nu|_{\hat Z}$.

Recall that $\psi := \frac{d\nu}{d\mu_\delta}$. By projecting
$\hat\nu_{\hat Z}$ down to the Julia set, we find that $\psi \ge
\hat \psi \circ \pi^{-1} > 0$ on $Z$. Let $\psi_0 = \inf\{ \psi(z)
: z \in Z\}$. Since $Z = U \cap \J$ for some open set $U$ in $\C$,
we can find $M$ \st $f^M(U) \supset \J$. Let us now prove that
$\psi > 0$ for other points as well, let $z \in \J \setminus
\cup_{i=1}^M f^i(\Crit)$ be arbitrary, and let $B \owns z$ be a
neighbourhood of $z$ such that $\diam(B) \le d(B,\cup_{i=1}^M
f^i(\Crit))$. There there is a subset $B_0 \subset U$ such that
$f^M:B_0 \to B$ is univalent. It follows that
\[
\nu(B) \ge \nu(B_0) \ge \psi_0 \mu_\delta(B_0) \ge \psi_0 \inf\{
|Df^M(z)|^{-\delta} : z \in Z \} \mu_\delta(B).
\]
This implies that $\psi(z) \ge  \psi_0 \inf\{ |Df^M(z)|^{-\delta}
: z \in Z \} > 0$.
\end{proof}

The following result clarifies some properties of $\delta$ and
$\mu_\delta$. The uniqueness part  is due to \cite{DMNU} and parts
(a) and (b) are due to \cite[Lemma 4.2]{BMO}. We let
\[
L(f):= \bigcup_{M>0}\bigcup_{\delta>0} \{z\in \J \hbox{ goes to
large scale for } \delta,M>0 \mbox{ infinitely often}\}.
\]
Points in this set are often referred to as \emph{conical} points.
For a system $(X,T,\mu)$ we say that $A$ is \emph{lim sup full} if
$\limsup_n\mu(T^nA)=1$.  We say that $T$ is lim sup full if this
property holds for all sets of positive measure.

\begin{theorem}
Suppose that $\mu_\delta$ is a $\delta$--conformal measure with
$\mu_\delta(L(f))>0$.  The $\mu_\delta$ is the unique measure with
this property and
\begin{itemize}
\item[(a)] $f$ is lim sup full, exact, ergodic,
conservative, $\mu$ is non--atomic, $\supp (\mu)=\J$ and
$\omega(z) =\J$ for $\mu$--a.e. $z \in \C$;
\item[(b)] $\delta$ is the minimal exponent for which a conformal
measure with support on $\J$ exists.
\end{itemize}
\label{thm:bmo}
\end{theorem}

Note that (b) implies the well-known fact that for any $f$ and
$\delta > 0$ satisfying the conditions of the theorem,
$\mu_{\delta'}(L(f))=0$ for each $\delta'>\delta$.
Mayer \cite{May} gives an example of a polynomial $f$ such that
$\mu_\delta(L(f)) = 0$ for all $\delta$ such that
$\delta$-conformal measure $\mu_\delta$ exists.

We next make an alternative assumption on the behaviour of points
under iteration by $f$ which guarantees that there is some lifted
measure.

\begin{theorem}
Suppose that \eqref{CO} is satisfied.  Let $\mu_\delta$ be a
$\delta$--conformal measure on $\J$, then the following are
equivalent.
\begin{itemize}
\item[(a)] There exists $\lambda>0$ \st
$\ulambda(z)\ge\lambda$ for all $z$ in a set of positive
$\mu_\delta$--measure;
\item[(b)] The measure $\mu_\delta$ is liftable.
\end{itemize}
\label{thm:pointwise lyap}
\end{theorem}

All of the situations considered in \cite{GS1,Prz,Rees} give
invariant probability measures $\mu \ll \mu_\delta$ for some
$\delta$--conformal measure with $\lambda_f(\mu)>0$.  Therefore,
all of those cases fit into our setting.  The closest result to
ours that we know of is \cite{GS1} where the measure $\mu$ was
obtained whenever the rational function $f$ satisfied a
summability condition on the derivatives of critical orbits.

\begin{corollary}
If there is a liftable probability measure $\mu \ll \mu_\delta$,
then $\delta = \dim_H(\mu)$ (where $\dim_H$ stands for Hausdorff
dimension).
\end{corollary}

\begin{proof}
Since $\mu$ is liftable, so is $\mu_\delta$.  By
Lemma~\ref{lem:conformal equiv}, $\mu$ is ergodic and by
Theorem~\ref{thm:pointwise lyap}, $\mu$ must have positive
Lyapunov exponent.  Pesin's formula in \cite[Chapter 10]{PU}
implies that $\delta = \dim_H(\mu)$.
\end{proof}

To prove Theorem~\ref{thm:pointwise lyap}, we will need the
following results. Define $s_R(n,D)$ to be the maximal number of
$n$--paths originating from an element $D \in \D,\ \lev (D)=R$ and
not re--entering $\hat\J_R$.  Let
\[
s_R(n):= \max \{s_R(n, D) : \lev(D) = R \}.
\]

\begin{lemma}\label{lem:rnbound}
Let $N := \# \P_1 \le \sum_c \kappa_c d_c$. For each
$R$, there exists a $C>0$ \st for all $0 \le j < R$,
$$s_R( nR+j ) \le C(2RN)^{n+1}.$$
\end{lemma}

The proof of this lemma is in the appendix.

\begin{proof}[Proof of Theorem~\ref{thm:pointwise lyap} assuming
Lemma~\ref{lem:rnbound}.] First assume that (b) holds, and let
$\hat \mu_\delta$ be the lifted measure.  By
Lemma~\ref{lem:invmeas}, $\hat\mu_\delta$ is invariant and by
Lemma~\ref{lem:conf_ac}, $\hat\mu_\delta\circ\pi^{-1}\ll\mu$.
Therefore Proposition~\ref{prop:lyap equal} implies that (a)
holds.

Now assume that (a) holds.  We will use a counting argument to
prove that a positive measure set of points must return to some
$\hat \J_R$ with positive frequency, from which liftability
follows.

For $1< \lambda_0 < \lambda$, $R, n \ge 1$ and $\eps>0$ we
consider the set
\[
B_{\lambda_0,R,n}(\eps):=\left\{z: |Df^n(z)|>\lambda_0^n \hbox{
and } \frac 1n\#\left\{0\le j<n:\hat f^j(i(z)) \in
\hat\J_R\right\}\le\eps\right\}.
\]
We let $\P_{B,n}$ denote the collection of cylinder sets of $\P_n$
which intersect $B_{\lambda_0,R,n}$. Since $\mu_\delta$ is
$\delta$--conformal, we can compute that
$\mu_\delta(B_{\lambda_0,R,n}(\eps)) \le \lambda_0^{-\delta
n}\#\P_{B,n}$. We will prove that by taking $R_0 \ge 1$ and
$\eps>0$ appropriately, this is arbitrarily small in $n$, which
leads us to conclude that a positive measure set must visit
$\hat\J_{R_0}$ with positive frequency.

Notice that any $Z \in \P_n$ uniquely determines a path $D_0(\hat
Z) \to \cdots \to D_{n-1}(\hat Z)$ in $\hat \J$ given by $\hat Z =
i(Z)$, $\hat f^j(\hat Z) \subset D_j(\hat Z) \in \D$, and vice
versa. In our case, given $P \in \P_{B,n}$, we let $\hat P =
i(P)$, we have a path defined in $\hat \J$. Moreover, $D_j(\hat P)
\cap \hat\J_R = \es$ for at most $\eps n$ of the times $j=0,
\ldots, n-1$.  We will estimate $\#\P_{B,n}$ in terms of these
paths. Define
$$
S(\eps, n) := \left\{M \subset \{0, \ldots , n-1\}: \# M \le \eps
n\right\}.
$$
The following well-known result estimates the cardinality of this
set . For $x\in (0,1)$, define $l(x):= -x\log x - (1-x)\log
(1-x)$.

\begin{lemma}\label{lem:stirling}
Let $S(\eps, n) := \left\{ M \subset \{0, \ldots , n-1\}: \# M \le
\eps n\right\}$.  Then for $n$ large, $\# S(\eps, n) \le
e^{n(\eps+l(\eps))}$.
\end{lemma}

Observe that for $M \in S(\eps, n)$, the set $\{0, \ldots, n-1\}
\setminus M$ consists of at most $1+\#M$ integer--intervals. The
number of $1$--paths in $\hat \J_R$ is bounded by the number of
domains in $\hat \J_{R+1}$.  By Lemma~\ref{lem:basiccounting}(c),
this is bounded above by $1+(R+1)\#\Crit\prod_c\kappa_c$.  Then choosing some
large $n \ge 1$,
\begin{align*}
\#\P_{B,n} &\le \sum_{M \in S(\eps, n)} \# \{Z \in \P_n: j \notin
M \Rightarrow D_j(i(Z))\cap \hat \J_R
= \es \} \\
& \le \sum_{M \in S(\eps, n)} \left(1+(R+1)\#\Crit\prod_c\kappa_c\right)^{\#M}
s_R(n-\#M) \\
& \le \#S(\eps, n) \left(1+(R+1)\#\Crit\prod_c\kappa_c\right)^{\eps n} s_R(n).
\end{align*}
Therefore, using Lemma~\ref{lem:rnbound}, there exist $R_0\ge 1$
and $\eps_0>0$ \st for some $\gamma<\lambda_0^\delta$ and $C>0$,
$\#\P_{B,n}< C\gamma^n.$  Therefore we have
$\mu_\delta\left(B_{\lambda_0,R_0,n}(\eps_0)\right) \le
C\left(\frac{\gamma}{\lambda_0^\delta}\right)^n$.  Whence,
$\mu_\delta\left(B_{\lambda_0,R_0,n}(\eps_0)\right) \to 0$ as $n
\to \infty$.  Since, by assumption, we have $\lim_{n \to
\infty}\mu_\delta\left\{z: |Df^n(z)|>\lambda_0^n\right\} > 0$,
there must exist some $\eps_1,\alpha>0$ \st for large enough $n\ge
1$,
\[
\mu_\delta\left\{z: \frac 1n\#\left\{0\le j<n:\hat f^j(i(z)) \in
\hat\J_{R_0}\right\}>\eps_1\right\}>\alpha.
\]
It is now easy to see that for any vague limit $\hat\mu_\delta$ of
measures obtained as in \eqref{eqn:mulift}, we have
$\hat\mu_\delta (\hat \J)\ge \hat\mu_\delta (\hat
\J_{R_0})>\alpha\eps_1$.
\end{proof}

\section*{Appendix}

This appendix is devoted to proving Lemma~\ref{lem:rnbound}. We
fix some $R \ge 1$ and a domain $D \in \D,\ \lev(D)=R$. We say
that
\begin{itemize}
\item a $t$--path {\em survives} if the path starts in $D$ and
never falls into $\hat \J_R$;
\item a domain {\em is surviving at time $t$} if it is the
terminal domain of a surviving $t$--path;
\item a cutpoint $z$ is a {\em surviving cutpoint at time $t$} if it
lies in the terminal domain of a surviving $t$--path.
\end{itemize}

Define
\begin{equation*}
L_t(m) := \#\{ \mbox{surviving $m$--cutpoints in $t$--paths
starting from } D\}.
\end{equation*}
Since for a path $D\to \cdots \to D'$, each cutpoint in $D$ has
only one image in $D'$, we have
\begin{equation}\label{eqn:lookback}
L_t(m) \le L_{t-l}(m-l) \mbox{ for } 1 \le l \le m \le t,
\end{equation}
which is a rule we will apply repeatedly. Moreover, since the
terminal domain of each surviving $t$--path contains at least one
$l$--cutpoint for $R < l \le R+t$, we find $L_t(j) = 0$ for $j >
t+R$ and
\begin{equation}\label{eqn:1cutpoints}
L_t(1) \le N \cdot \#\{\mbox{surviving $t-1$--paths}\} \le N
\left(\sum_{l=R+1}^{R+t-1} L_{t-1}(l)\right),
\end{equation}
where $N=\# \P_1$.
Using these rules, we prove the following lemma.
\begin{lemma}\label{lem:findrnbound}
Suppose that $(n-1)R < t \le nR$, then
\begin{equation}\label{bound_Lt}
L_t(j) \le \left\{ \begin{array}{ll}
2^n R^{n} N^{n+1} & \mbox{ if } 0 < j \le t;\\
N  & \mbox{ if } t < j \le t+R; \\
0  & \mbox{ if } t+R < j.
\end{array} \right.
\end{equation}
\end{lemma}

\begin{proof}
Since every terminal domain of an $t$--path has level $\le t+R$,
$L_t(j) = 0$ for $j > t+R$. This proves the third inequality.
Before we prove the remaining part by induction, let us
compute what happens for $t \le R$. \\[0.2cm]
\underline{\bf $t=0$:} The maximal number of $1$--cutpoints
possible in a single domain is $\#\Crit \le N$. So in particular
$L_0(1) \le N$. By Lemma~\ref{lem:basiccounting}, $L_0(l) \le N$
for
$1 \le l \le R$.\\[0.2cm]
\underline{\bf $t=1$:} By rule \eqref{eqn:lookback}, $L_1(j+1) =
L_0(j) \le N$ for $1 \le j \le R$. By rule
\eqref{eqn:1cutpoints}, $L_1(1) \le N$.\\[0.2cm]
\underline{\bf $t=2$:} By rule \eqref{eqn:lookback}, $L_2(2+j) =
L_0(j) \le N$ for $1 \le j \le R$. Similarly, $L_2(2) = L_1(1) \le
N$.  Also $L_2(1) \le N^2$ by rule \eqref{eqn:1cutpoints}.
\\[0.2cm]
\underline{\bf $t \le R$:} As before
\[
L_t(j) = L_0(j-t) \le N \mbox{ for } t < j \le t+R,
\]
and
\begin{align*}
L_t(j) &\le L_{t-j+1}(1) \le N \# \{ \mbox{$t-j$--surviving paths} \}\\
&\le N \sum_{l=R+1}^{t-j+R} L_{t-j}(l) \le (t-j)N^2 \le RN^2
\qquad \mbox{by rule \eqref{eqn:1cutpoints},}
\end{align*}
for $0 < j \le t$.\\[0.2cm]
It follows that if $t < j \le t+R$, then $L_t(j) = L_0(j-t) \le
N$. So now we have proved the second inequality of the lemma for
all $t$, and the first inequality for $t \le R$.

We continue by induction on $n$. So assume that \eqref{bound_Lt}
holds for $n$ and that $nR < t \le (n+1)R$ and $j \le t$. Then
\begin{align*}
L_t(j) &= L_{t-j+1}(1) \le N\sum_{l = R+1}^{t-j+R} L_{t-j}(l)
\qquad \qquad \mbox{ by rule \eqref{eqn:1cutpoints}}
\\
&\le RN^2+ N \sum_{l = R+1}^{t-j} L_{t-j}(l)\hspace{2.5mm}
\mbox{by the induction hypothesis for $l > t-j$}\\
\qquad &\le  RN^2+ N \sum_{l = R+1}^{t-j} L_{t-j-l+1}(1)
\qquad \mbox{ by rule \eqref{eqn:lookback}}\\
&\le RN^2+ N \sum_{s = 1}^{t-j-R} L_s(1) \\
&\le RN^2+ N \sum_{d=1}^n \sum_{s = (d-1)R+1}^{dR} L_s(1)\\
&\le RN^2+ N \sum_{d=1}^n \sum_{s = (d-1)R+1}^{dR} 2^d R^{d}
N^{d+1} \hspace{3mm} \mbox{by the induction hypothesis}\\
&\le RN^2+ RN^2 \sum_{d=1}^n (2RN)^d \le RN^2
\left(1+ \frac{(2RN)^{n+1}-2RN}{2RN-1}\right) \\
& \le  RN^2\left(\frac{(2RN)^n}{1-\frac{1}{2RN}}\right) \le
2^{n+1} R^{n+1} N^{n+2}.
\end{align*}
This proves the induction step, and hence the lemma.
\end{proof}

\begin{proof}[Proof of Lemma~\ref{lem:rnbound}.]
Since there cannot be more surviving $nR$--paths than surviving
cutpoints at time $nR$ we can estimate $s_R$ using $L_t(j)$. We
use rules \eqref{eqn:lookback} and \eqref{eqn:1cutpoints} as in
the previous proof.
\begin{align*}
s_R( nR+j ) &= \#\{ \mbox{$nR+j$--surviving paths} \} \\
&\le \sum_{l = R+1}^{nR+j+R} L_{nR+j}(l) \qquad  \qquad
\mbox{by rule \eqref{eqn:1cutpoints}}\\
&\le RN^2+ \sum_{l = R+1}^{nR+j} L_{nR+j-l+1}(1)
\qquad \mbox{by rule \eqref{eqn:lookback} and Lemma~\ref{lem:findrnbound}}\\
&\le RN^2+ \sum_{s = 1}^{(n-1)R+j} L_s(1) \\
&\le RN^2+ \sum_{d=1}^n \sum_{s = (d-1)R+1}^{dR} L_s(1) \\
&\le RN^2+ \sum_{d=1}^n \sum_{s = (d-1)R+1}^{dR} (2RN)^{d+1}
\qquad \mbox{by Lemma~\ref{lem:findrnbound}}\\
&\le C (2RN)^{n+1},
\end{align*}
for some constant $C > 0$.
\end{proof}

\medskip
\noindent
Department of Mathematics\\
University of Surrey\\
Guildford, Surrey, GU2 7XH\\
UK\\
\texttt{h.bruin@surrey.ac.uk}\\
\texttt{http://www.maths.surrey.ac.uk/showstaff?H.Bruin}

\medskip
\noindent
Department of Mathematics\\
University of Surrey\\
Guildford, Surrey, GU2 7XH\\
UK\\
\texttt{m.todd@surrey.ac.uk}\\
\texttt{http://www.maths.surrey.ac.uk/showstaff?M.Todd}

\medskip

\end{document}